\documentclass[12pt]{article}

\usepackage{amssymb,amsmath,amsthm,framed}

\newcommand{\ovl}{\overline}
\newcommand{\n}{\noindent}

\newcommand{\cl}[1]{{\mathcal{#1}}}
\newcommand{\bb}[1]{{\mathbb{#1}}}
\newcommand{\ovltimes}{\ovl\otimes}

\numberwithin{equation}{section}

\theoremstyle{plain}
\newtheorem{lem}{Lemma}[section]

\newtheorem{thm}[lem]{Theorem}
\newtheorem{cor}[lem]{Corollary}

\theoremstyle{definition}

\newtheorem{exm}[lem]{Example}

\newtheorem{rem}[lem]{Remark}

\begin{document}

\title{NORMALIZERS OF IRREDUCIBLE SUBFACTORS}

\author{Roger R. Smith\thanks{Partially supported by NSF grant DMS-0401043.}, \and Stuart A. White, \and Alan D. Wiggins\\
\phantom{blankline}\\
\phantom{blankline}\\
Department of Mathematics\\
Texas A\&M University\\
College Station, TX 77843\\
\phantom{blankline}\\
rsmith@math.tamu.edu, white@math.tamu.edu, awiggins@math.tamu.edu}

\date{}
\maketitle
\begin{abstract}
We consider normalizers of an irreducible inclusion $N\subseteq M$ of $\mathrm{II}_1$ factors.   In the infinite index setting an inclusion $uNu^*\subseteq N$ can be strict, forcing us to also investigate the semigroup of one-sided normalizers.  We relate these normalizers of $N$ in $M$ to projections in the basic construction and show that every trace one projection in the relative commutant $N'\cap \langle M,e_N\rangle$ is of the form $u^*e_Nu$ for some unitary $u\in M$ with $uNu^*\subseteq N$. This enables us to identify the normalizers and the algebras they generate in several situations. In particular each normalizer of a tensor product of irreducible subfactors is a tensor product of normalizers modulo a unitary. We also examine normalizers of irreducible subfactors arising from subgroup--group inclusions $H\subseteq G$. Here the normalizers are the normalizing group elements modulo a unitary from $L(H)$. We are also able to identify the finite trace $L(H)$-bimodules in $\ell^2(G)$ as double cosets which are also finite unions of left cosets. 
\end{abstract}

\thispagestyle{empty}

\newpage

\setcounter{page}{1}

\section{Introduction}\label{SWWsec1}

\indent 

In this paper we consider the following general problem:\ given an irreducible inclusion $N\subseteq M$ of $\text{II}_1$ factors, is it possible to determine the normalizing unitaries and the von Neumann algebra that they generate? In such generality very little can be said, but we will examine two naturally occurring classes of inclusions where this question can be answered completely. The first of these arises from an inclusion $H\subseteq G$ of countable discrete groups. Certain algebraic conditions ensure that the inclusion $L(H)\subseteq L(G)$ of group von Neumann algebras is irreducible, and when $H$ has finite index in $G$ we have one of the most basic examples in subfactor theory. We will not assume finite index, and so our results also apply to the infinite index case. The second situation requires a finite group $G$ to act by outer automorphisms on a $\text{II}_1$ factor $M$. The fixed point algebra $M^G$ gives a finite index irreducible inclusion $M^G\subseteq M$ for which the algebra generated by the normalizing unitaries can be determined. The  lemmas that we develop for the solution of this problem also allow us to describe the normalizers of tensor products in terms of tensor products of normalizers. The main technical result is Theorem \ref{SWWthm3.4} which relates one-sided normalizers to projections in $N'\cap \langle M,e_N\rangle$ of unit trace. This is the basis for the main theorems of the paper, which are Theorems \ref{SWWthm4.2}, \ref{SWWthm5.4}, \ref{SWWthm6.2} and \ref{SWWthm7.1}. We describe these in more detail below. 

Dixmier, \cite{Dixmier.Masa}, was the first to recognize the importance of the normalizer ${\cl N}(A)$ for a von Neumann subalgebra $A$ of a factor $M$. In the case of maximal abelian self-adjoint subalgebras (masas), he classified the masas according to whether ${\cl N}(A)''$ was $M$ (regular), was a proper subfactor (semiregular), or was equal to $A$ (singular). He also provided examples of each type by considering inclusions $H\subseteq G$ of suitably chosen group-subgroup pairs. Masas satisfy $A=A'\cap M$ and so their commutants are large. The opposite end of the spectrum is the condition $N'\cap M={\bb C}1$, which defines an irreducible subfactor. Such subfactors will be the focus of our study. The isolated examples of singular subfactors in \cite{Sinclair.StrongSing2} were the starting point for a systematic examination of this phenomenon for inclusions of the form $M\rtimes_\alpha H \subseteq M\rtimes_\alpha G$ in \cite{Saw.Length,Wiggins.StrongSing}. Singularity was connected to strong singularity of masas in \cite{Sinclair.StrongSing,Saw.StrongSing}, and one consequence of this was the formula
\begin{equation}\label{SWWeq1.1}
{\cl N}(A_1\ \ovl\otimes\ A_2)'' = {\cl N}(A_1)''\ \ovl\otimes\ {\cl N}(A_2)''
\end{equation}
of \cite{Saw.StrongSing} for singular masas, which simply says that the tensor product of singular masas is again singular. Subsequently Chifan, \cite{Chifan.Normalisers}, proved (\ref{SWWeq1.1}) for general masas. These papers collectively have provided strong motivation for the work undertaken here. It also depends heavily on the recent theory of perturbations, developed primarily by Popa, \cite{Popa.Betti,Popa.StrongRigidity1,Popa.StrongRigidity2,Sinclair.PertSubalg,Popa.IPP}, building of the work of Christensen \cite{Christensen.Subalgebras} in which an important averaging technique is developed.  The normalizers of regular masas play a key role in this perturbation theory. A second crucial ingredient is the theory of subfactors \cite{Jones.Index,Jones.SubfactorsBook,Popa.CBMSNotes,Popa.Entropy}.

The link to subfactor theory is made by the simple observation that if $u\in {\cl N}(N)$, then the operator $u^*e_Nu$ lies in the relative commutant $N'\cap \langle M,e_N\rangle$ for the basic construction $\langle M,e_N\rangle$ (see the next section for explanations of terminology). It is then appropriate to study normalizers in the context of the basic construction and this relative commutant, which has proved to be a fruitful interaction. There is also a further ingredient:\ an unexpected role in our work is played by the more general one-sided normalizers, those unitaries $u\in M$ satisfying $uNu^*\subseteq N$. This containment can be strict, so normalizers and their one-sided counterparts are distinct in general, as we show in Example \ref{SWWexm6.3}.  

For a masa $A\subseteq M$, any unitary $u$ which is a one-sided normalizer of $A$ has the property that $uAu^* \subseteq uMu^* = M$ is a masa in $M$ contained inside the masa $A$. The defining property of masas then implies that $uAu^*=A$ and $u$ is also a normalizing unitary. A second situation where equality occurs is a finite index inclusion of factors $N\subseteq M$. Each one-sided normalizing unitary $u$ induces an equivalence of containments $N\subseteq M$ and $uNu^* \subseteq uMu^* = M$ which then have equal finite indexes. This is incompatible with $uNu^*\subseteq N\subseteq M$ unless the first two algebras are equal, in which case $u$ is a normalizing unitary. As a consequence, the one-sided normalizers do not appear until we consider infinite index inclusions of factors, and these form a substantial part of the paper.

The contents of the paper are as follows. Section \ref{SWWsec2} establishes notation, reviews some well known facts about the basic construction, and recalls the Galois theory for finite group actions. Section \ref{SWWsec3} examines the interplay between one-sided normalizers and projections in $N'\cap \langle M,e_N\rangle$ when $N$ is irreducible. Here it is shown that every such projection $f$ satisfies $\text{Tr}(f)\ge 1$, and is of the form $u^*e_Nu$ for a one-sided normalizer $u$ precisely when $\text{Tr}(f)=1$. These results occur in Theorem \ref{SWWthm3.4} which is the technical basis for our subsequent theorems. Section \ref{SWWsec4} characterizes both one-sided normalizers and normalizers for tensor products (Theorems \ref{SWWthm4.1} and \ref{SWWthm4.2}), while Section \ref{SWWsec5} considers normalizers for the inclusion $M^G\subseteq M$, where $G$ is a finite group. In this case we have a finite index inclusion, and so there is no distinction between normalizers and one-sided normalizers, as noted above. The algebra generated by the normalizers of $M^G$ is shown to be $M^{[G,G]}$ where $[G,G]$ is the commutator subgroup. Surprisingly, it has been possible to describe this algebra without being able to identify individual non-trivial normalizing unitaries in this case.

The last two sections are devoted to group-subgroup inclusions. When $L(H) \subseteq L(G)$, we characterize the normalizers and one-sided normalizers of $L(H)$ in terms of their counterparts at the group level. The ranges of projections in $L(H)' \cap \langle L(G), e_{L(H)}\rangle$ are the $L(H)$-bimodules in $\ell^2(G)$. We investigate these in Section \ref{SWWsec7}. Those that correspond to projections of finite trace are characterized algebraically in terms of left cosets and double cosets in Theorem \ref{SWWthm7.1}, while the subsequent examples  show that the situation is much more complicated for projections of infinite trace.

The following useful analogy between masas and subfactors has been implicit in much of the last two sections. For a masa $A\subseteq M$, the Puk\'anszky invariant is defined by using the algebra ${\cl A}' = (A\cup JAJ)'$, and this can also be viewed as the relative commutant $A'\cap \langle M,e_A\rangle$. It is type ${\rm{I}}$, and the integers (including $\infty$) which comprise the Puk\'anszky invariant come from the various summands of type ${\rm{I}}_n$ in $e^\bot_A{\cl A}'$. For irreducible inclusions of factors $N\subseteq M$, essentially the same algebra $N'\cap \langle M,e_N\rangle$ occurs, where $e_N$ is central just as $e_A$ is central in the masa case. When an abelian subgroup $H\subseteq G$ generates a masa $L(H)$ in $L(G)$, it is often the case that the Puk\'anszky invariant can be determined from the structure of the double cosets $HgH$ in $G$ \cite{Sinclair.Puk,Sinclair.FreePuk}. These may be identified with $L(H)$-bimodules in $\ell^2(G)$, and as such they play a significant part in Section \ref{SWWsec7} where subfactors arising from subgroups are considered. The interplay between these various quantities has been studied extensively in the theory of finite index inclusions of factors \cite{Jones.Index,Jones.SubfactorsBook,Kosaki.Bimodules} but the methods developed there do not seem helpful for the infinite index situation.

We end by pointing out that Theorem \ref{SWWthm6.2}  could have been deduced from Theorem \ref{SWWthm7.1}. We have chosen not to do so, because the techniques for Theorem \ref{SWWthm7.1} are more complicated and appear in a more understandable form in the proof of Theorem \ref{SWWthm6.2}.

\section{Notation and preliminaries}\label{SWWsec2}

\indent 

The basic object of study in this paper is an inclusion $N\subseteq M$ of $\text{II}_1$ factors, where the unique normalized faithful normal trace on $M$ is denoted by $\tau$. We will always assume that these factors are separable although this is just for notational convenience; the results are valid in general. We always assume that $M$ is in standard form, so that it is represented as left multiplication operators on the Hilbert space $L^2(M,\tau)$, or simply $L^2(M)$. We reserve the letter $\xi$ to denote the image of $1\in M$ in this Hilbert space, and $J$ will denote the isometric conjugate linear operator on $L^2(M)$ defined by
\begin{equation}\label{SWWeq2.1}
J(x\xi) = x^*\xi,\qquad x\in M,
\end{equation}
and extended by continuity to $L^2(M)$ from the dense subspace $M\xi$. Then $L^2(N)$ is a closed subspace of $L^2(M)$, and $e_N$ denotes the projection of $L^2(M)$ onto $L^2(N)$. The basic construction is the von Neumann algebra generated by $M$ and $e_N$, and is denoted $\langle M,e_N\rangle$. Since $M'\cap B(L^2(M)) = JMJ$, we also have $\langle M,e_N\rangle' \cap B(L^2(M)) = JNJ$. This shows that $\langle M,e_N\rangle$ is either type $\text{II}_1$ or $\text{II}_\infty$, and in both cases there is a unique semifinite normal trace Tr with the property that $\text{Tr}(e_N)=1$. The Jones index can be described as $\text{Tr}(1)$, although this is not the original definition. These are standard facts in subfactor theory, and can be found in \cite{Jones.Index,Popa.Entropy,Jones.SubfactorsBook}. These sources also contain the following properties of the Jones projection $e_N$ which we now list. We will use them subsequently without comment. The unique trace preserving conditional expectation of $M$ onto $N$ is denoted ${\bb E}_N$.
\begin{itemize}
\item[(i)] $e_N(x\xi) = {\bb E}_N(x)\xi$, \ $x\in M$.
\item[(ii)] $e_Nxe_N = {\bb E}_N(x)e_N = e_N{\bb E}_N(x), \ x\in M$.
\item[(iii)] $x\mapsto e_Nx$ and $x\mapsto xe_N$ are injective maps for $x\in M$.
\item[(iv)] $\{xe_Ny\colon \ x,y\in M\}$ generates a strongly dense subalgebra of $\langle M,e_N\rangle$.
\item[(v)] $\text{Tr}(xe_Ny) = \tau(xy)$ for $x,y\in M$.
\item[(vi)] $Me_N$ is $*$-strongly dense in $\langle M,e_N\rangle e_N$.
\item[(vii)] $e_N\langle M,e_N\rangle e_N=Ne_N=e_NN$.
\item[(viii)] $M\cap\{e_N\}'=N$.
\end{itemize}

One further fact which we will need is provided in the following lemma. It is surely known to experts, but we do not know a reference.

\begin{lem}\label{SWWlem2.1}
Let $N_i\subseteq M_i$, $i=1,2$, be inclusions of ${\rm II}_1$ factors
and let ${\text{Tr}}_i$ be the canonical trace on $\langle
M_i,e_{N_i}\rangle$, $i=1,2$.
 Then
\begin{equation}\label{SWWeq2.2}
\langle M_1,e_{N_1}\rangle\ \ovl\otimes\ \langle M_2,e_{N_2}\rangle \cong \langle M_1\ \ovl\otimes\ M_2, e_{N_1\ovl\otimes N_2}\rangle, 
\end{equation}
and ${\text{Tr}}_1\otimes {\text{Tr}}_2$ is the canonical
trace on the tensor product.
\end{lem}

\begin{proof}
We identify the Hilbert spaces $L^2(M_1)\otimes L^2(M_2)$ and $L^2(M_1\otimes M_2)$ and assume that both algebras in \eqref{SWWeq2.2} act on the former. This identifies
$e_{N_1\ovl\otimes N_2}$ with $ e_{N_1} \otimes e_{N_2}$. Tomita's commutant theorem gives that the commutants of these algebras are respectively
\[
(JN_1J)\ \ovl\otimes\ (JN_2J)\quad \text{and}\quad \left((JM_1J)\ \ovl\otimes\  (JM_2J)\right)\cap \{e_{N_1\ovl\otimes N_2}\}'.
\]
To show equality it suffices to apply $J\otimes J$ on the left and the right and instead prove equality of 
\[
N_1\ \ovl\otimes\ N_2\quad \text{and}\quad \left(M_1\ \ovl\otimes\  M_2\right)\cap \{e_{N_1\ovl\otimes N_2}\}'.
\]
But this is just (viii) above applied to the containment $N_1\ \ovl\otimes\ N_2\subseteq
M_1\ \ovl\otimes\ M_2$.
Because $({\text{Tr}}_1\otimes {\text{Tr}}_2)(e_{N_1}\otimes
e_{N_2})=1$, the statement on traces is now clear.
\end{proof}

When a finite group $G$ acts on a $\text{II}_1$ factor $M$ by outer automorphisms, the fixed point algebra $M^G$ is a subfactor of $M$. The Galois theory of \cite{Nakamura.GaloisFinite,Nakamura.Galois} then characterizes the intermediate subfactors in terms of subgroups in the following way:\ there is an order reversing one to one correspondence between subgroups $K$ of $G$ and intermediate subfactors given by $K\mapsto M^K$. We will use this in Section \ref{SWWsec5}.

The last two sections are concerned with inclusions $H\subseteq G$ of groups. The canonical basis for $\ell^2(G)$ is denoted $\{\delta_g\colon \ g\in G\}$, and we assume that $G$ is represented on this Hilbert space by the left regular representation $\lambda$, so that $\lambda_s\delta_t = \delta_{st}$ for $s,t\in G$. The right regular representation $\rho$ satisfies $\rho(s)=J\lambda_sJ$. As is standard, $L(G)$ is used for the von Neumann algebra generated by the left regular representation.

On several occasions we will require the following lemma. A similar result can be found in \cite{Saw.Length}, but this is not quite in the form that we need, so we offer a slightly more general version here.

\begin{lem}\label{SWWlem2.2}
Let $N\subseteq M$ be an inclusion of ${\rm{II}}_1$ factors on $L^2(M)$ such that $N'\cap M = {\bb C}1$. Let $\{\phi_1,\ldots,\phi_n\}$ be a set of automorphisms of $M$ with the property that the restriction of each $\phi^{-1}_j\phi_i$ to $N$ is not implemented by a unitary in $M$ whenever $i\ne j$. Let $X\subseteq N$ and $Y\subseteq JMJ$ be self-adjoint subsets which generate their respective containing factors and assume that $1\in X$. Then the von Neumann subalgebra of $\mathbb M_n(B(L^2(M)))$ generated by
\[
\left\{\begin{pmatrix}
\phi_1(x)y\\ &\ddots\\ &&\phi_n(x)y
       \end{pmatrix}\colon \ x\in X, y\in Y\right\}
\]
is
\[
\left\{\begin{pmatrix}
t_1\\ &\ddots\\ &&t_n
       \end{pmatrix}\colon \ t_i\in B(L^2(M))\right\}.
\]
\end{lem}

\begin{proof}
By the double commutant theorem, it suffices to show that the commutant of the first set of operators is the set of diagonal scalar matrices. Commutation with
\[
\left\{\begin{pmatrix}
y\\ &\ddots\\ &&y
       \end{pmatrix}\colon \ y\in Y\right\}
\]
allows us to consider a matrix $(m_{ij})\in{\bb M}_n(M)$. The conditions for this to commute with
\[
\left\{\begin{pmatrix}
\phi_1(x)\\ &\ddots\\ &&\phi_n(x)
       \end{pmatrix}\colon \ x\in X\right\}
\]
are
\begin{equation}\label{SWWeq2.3}
\phi_i(x) m_{ij} = m_{ij}\phi_j(x),\qquad x\in X,\quad 1\le i,j\le n,
\end{equation}
which then must hold for all $x\in N$. Since $\phi_i(N)$ has trivial relative commutant in $M$, while \eqref{SWWeq2.3} shows that $m_{ij}m^*_{ij} \in \phi_i(N)'\cap M$, we conclude that each $m_{ij}$ is a scalar multiple of a unitary. The case $i=j$ in \eqref{SWWeq2.3} places $m_{ii} \in \phi_i(N)'\cap M = {\bb C}1$ so the diagonal entries are scalars.

Suppose that some $m_{ij}\ne 0$ for a pair of integers $i\ne j$. By scaling we may replace $m_{ij}$ in \eqref{SWWeq2.3} by a unitary $u\in M$. If we  apply $\phi^{-1}_j$, then we obtain
\begin{equation}\label{SWWeq2.4}
\phi^{-1}_j(\phi_i(x)) \phi^{-1}_j(u) = \phi^{-1}_j(u)x,\qquad x\in N,
\end{equation}
which is contrary to the hypothesis that $\phi^{-1}_j\phi_i$ is not unitarily implemented on $N$. This shows that $m_{ij} = 0$ for $i\ne j$, completing the proof.
\end{proof}

\section{Normalizers and the basic construction}\label{SWWsec3}

\indent 

Throughout this section $N\subseteq M$ will denote an irreducible inclusion of $\text{II}_1$ factors. We make no assumption of finite index; when we restrict to the finite index situation, we will indicate this explicitly. Our aim is to relate the normalizers and one-sided normalizers of $N$ to the set of projections in the relative commutant $N'\cap \langle M,e_N\rangle$. 

The unitary group of $M$ is written as $\cl U(M)$ and we use the  notation   
\[
{\cl N}(N) = \{u\in {\cl U}(M)\colon \ uNu^*=N\},\ \  \cl{ON}(N) = \{u\in {\cl U}(M)\colon \ uNu^* \subseteq N\}
\]
to denote respectively the group of unitary normalizers and the semigroup of one-sided unitary normalizers of $N$.  The semigroup of one-sided normalizers is a sub-semigroup of the groupoid normalizers
$$
\cl{GN}(N)=\{v\in M:vNv^*\subseteq N,\quad vv^*,v^*v\in N\}.
$$
Since ${\cl N}(N) \subseteq \cl{ON}(N)$, we note that the following three lemmas are valid for normalizers in addition to one-sided normalizers.

\begin{lem}\label{SWWlem3.1}
Let $u\in \cl{ON}(N)$, and let $\phi\colon \ N\to N$ be the $*$-homomorphism defined by $\phi(x) = uxu^*$, $x\in N$. Then $e_N$ is a central projection in $\phi(N)' \cap \langle M,e_N\rangle$. In particular, this projection is central in $N'\cap \langle M,e_N\rangle$.
\end{lem}

\begin{proof}
Let $v$ be a fixed but arbitrary unitary in $\phi(N)' \cap \langle M,e_N\rangle$. We begin by establishing that $v\xi = \lambda\xi$ for some $\lambda\in {\bb C}$, $|\lambda|=1$.

Let $\eta = v\xi\in L^2(M)$. By $\|\cdot\|_2$-density of $M\xi$ in $L^2(M)$, we may find a sequence $\{x_n\}^\infty_{n=1}$ in $M$ such that
\begin{equation}\label{SWWeq3.1}
\lim_{n\to\infty}\|\eta-x_n\xi\|_2 = 0.
\end{equation}
Noting that $v$ commutes with $\phi(N)$ and with $JNJ$, we obtain
\begin{align}
J\phi(w) J\phi(w) \eta &= J\phi(w) J\phi(w) v\xi = vJ\phi(w) J\phi(w)\xi\notag\\
\label{SWWeq3.2}
&= v\xi = \eta,\qquad w\in {\cl U}(N),
\end{align}
so $\eta$ is an invariant vector for $J\phi(w)J\phi(w)$. Then
\begin{equation}\label{SWWeq3.3}
\|J\phi(w)J\phi(w) x_n\xi - \eta\|_2 = \|J\phi(w)J\phi(w)(x_n\xi-\eta)\|_2 \le \|x_n\xi-\eta\|_2.
\end{equation}
For each $n\in{\bb N}$, let $y_n\in M$ be such that $y_n\xi$ is the unique element of minimal $\|\cdot\|_2$-norm in
\[
\ovl{\rm conv}^w\{\phi(w)x_n\phi(w)^*\xi\colon \ w\in {\cl U}(N)\},
\]
by \cite[Section 8.2]{Sinclair.MasaBook}. Taking convex combinations and norm limits in \eqref{SWWeq3.3} shows that 
\begin{equation}
\|y_n\xi-\eta\|_2 \le \|x_n\xi-\eta\|_2,
\end{equation}
so \eqref{SWWeq3.1} implies that 
\begin{equation}\label{SWWeq3.4}
\lim\limits_{n\to\infty} \|y_n\xi-\eta\|_2 = 0.
\end{equation}
Moreover, uniqueness of $y_n\xi$ shows that
\begin{equation}\label{SWWeq3.5}
\phi(w) y_n\phi(w)^* =y_n,\qquad w\in {\cl U}(N),
\end{equation}
and so
\begin{equation}\label{SWWeq3.6}
\phi(w)y_n = y_n\phi(w),\qquad w\in {\cl U}(N).
\end{equation}
Thus $y_n\in \phi(N)'\cap M = {\bb C}1$ since $N$ is irreducible. From \eqref{SWWeq3.1} and \eqref{SWWeq3.4}, we conclude that $\eta =\lambda\xi$ for some $\lambda\in {\bb C}$. Since $\|\eta\|_2 = \|v\xi\|_2 =1$, it follows that $|\lambda|=1$.

For an arbitrary $x\in M$,
\begin{equation}\label{SWWeq3.7}
ve_Nx\xi = v{\bb E}_N(x)\xi = vJ{\bb E}_N(x^*)J\xi.
\end{equation}
Since $v$ commutes with $JNJ = \langle M,e_N\rangle'$, \eqref{SWWeq3.7} shows that
\begin{align}
ve_Nx\xi &= J{\bb E}_N(x^*)Jv\xi = J{\bb E}_N(x^*)J\lambda\xi\notag\\
\label{SWWeq3.8}
&= \lambda{\bb E}_N(x)\xi = \lambda e_Nx\xi,\qquad x\in M.
\end{align}
Thus $ve_N = \lambda e_N$, so
\begin{equation}\label{eqSWWeq3.9}
ve_Nv^* = (ve_N)(e_Nv^*) = \lambda e_N \bar\lambda e_N = e_N,
\end{equation}
since $|\lambda|^2 = 1$, showing that $v$ commutes with $e_N$. Since $v\in \phi(N)'\cap \langle M,e_N\rangle$ was an arbitrary unitary, we conclude that $e_N$ is central in this algebra.

The second statement of the lemma is an immediate consequence of taking $u$ to be 1, whereupon $\phi(N)=N$.
\end{proof}

\begin{lem}\label{SWWlem3.2}
Let $u\in \cl{ON}(N)$ be a fixed but arbitrary unitary. Then $u^*e_Nu$ is a minimal projection in $N'\cap \langle M,e_N\rangle$ and is also central in this algebra.
\end{lem}

\begin{proof}
As in the proof of Lemma \ref{SWWlem3.1}, let $\phi\colon \ N\to N$ be the $*$-homomorphism defined by $\phi(x)=uxu^*$ for $x\in N$. If $y\in N'\cap \langle M,e_N\rangle$ and $x\in N$,
\begin{equation}\label{SWWeq3.10}
(u^*yu)x = u^*y\phi(x)u = u^*\phi(x)yu = x(u^*yu).
\end{equation}
Then \eqref{SWWeq3.10} shows that $u^*yu\in N'\cap \langle M,e_N\rangle$ whenever $y\in N'\cap \langle M,e_N\rangle$. In particular, $u^*e_Nu \in N'\cap \langle M,e_N\rangle$ by letting $y$ be $e_N$.

To establish minimality of $u^*e_Nu$, consider a projection $q\in N'\cap \langle M,e_N\rangle$ satisfying $q\le u^*e_Nu$. Then $uqu^*\le e_N$, so there is a projection $p\in N$ such that $uqu^* = pe_N$, as implied by the relation $e_N\langle M,e_N\rangle e_N=Ne_N$ from Section \ref{SWWsec2}. For each $x\in N$,
\begin{align}
\phi(x) pe_N &= uxu^* uqu^* = uxqu^* = uqxu^*\notag\\
\label{SWWeq3.11}
&= uqu^*\phi(x) = pe_N\phi(x) = p\phi(x)e_N,
\end{align}
and so $\phi(x)p  = p\phi(x)$ for $x\in N$. Thus  $p\in \phi(N)'\cap N\subseteq \phi(N)'\cap M ={\bb C}1$, showing that $p=0$ or $p=1$. It follows that $q=0$ or $q=u^*e_Nu$, proving minimality of $u^*e_Nu$ in $N'\cap \langle M,e_N\rangle$.

We now show centrality of $u^*e_Nu$. For $z\in N$ and $y\in N'\cap\langle M,e_N\rangle$,
\begin{equation}\label{SWWeq3.12}
uyu^*\phi(z)=\phi(z)uyu^*
\end{equation}
by taking $x=\phi(z)$ in (\ref{SWWeq3.10}). So $uyu^*\in\phi(N)'\cap \langle M,e_N\rangle$. By Lemma \ref{SWWlem3.1}, $e_N$ is central in the latter algebra from which we obtain 
\begin{equation}\label{SWWeq3.13}
e_Nuyu^* = uyu^*e_N.
\end{equation}
In \eqref{SWWeq3.13}, multiply on the left by $u^*$ and on the right by $u$. The result is that $u^*e_Nu$ commutes with $y$ for $y\in N'\cap \langle M,e_N\rangle$, showing centrality of $u^*e_Nu$.
\end{proof}

The next lemma determines when $u^*e_Nu$ and $v^*e_Nv$ are equal for unitaries $u,v\in \cl{ON}(N)$.

\begin{lem}\label{SWWlem3.3}
Let $u,v\in \cl{ON}(N)$. Then $u^*e_Nu = v^*e_Nv$ if and only if there exists $w\in {\cl U}(N)$ such that $v=wu$.
\end{lem}

\begin{proof}
If $v=wu$ for some $u\in {\cl U}(N)$, then
\begin{equation}\label{SWWeq3.14}
v^*e_Nv = u^*w^*e_Nwu = u^*e_Nu,
\end{equation}
since $w$ commutes with $e_N$.

Conversely, suppose that $u^*e_Nu = v^*e_Nv$. Then $vu^*e_Nuv^* =e _N$. Multiply this equation on the right by $e_N$ to obtain $vu^*{\bb E}_N(uv^*) e_N=e_N$. Since $vu^*{\bb E}_N(uv^*)\in M$, we conclude that $vu^*{\bb E}_N(uv^*) = 1$, so $uv^*= {\bb E}_N(uv^*)$, and ${\bb E}_N(uv^*)$ is a unitary in $N$, which we denote by $w^*$. Then $v=wu$ as required.
\end{proof}

The following theorem is the main technical result which links one-sided normalizers to certain projections in $N'\cap \langle M,e_N\rangle$. The construction of the $*$-homomorphism $\phi$ in the proof comes from \cite{Popa.Betti}.

\begin{thm}\label{SWWthm3.4}
\begin{itemize}
\item[\rm (i)] Each non--zero projection of $f\in N'\cap \langle M,e_N\rangle$ satisfying ${\rm Tr}(f)\le 1$ has the form $u^*e_Nu$ for some $u\in \cl{ON}(N)$. 
\item[\rm (ii)] Each non--zero projection $f\in N'\cap \langle M,e_N\rangle$ satisfies ${\rm Tr}(f)\ge 1$.
\end{itemize}
\end{thm}

\begin{proof}
By Lemma \ref{SWWlem3.2}, each projection $u^*e_Nu$, where $u\in \cl{ON}(N)$, lies in $N'\cap \langle M,e_N\rangle$. Conversely, let $f$ be a non--zero projection in $N'\cap \langle M,e_N\rangle$ satisfying $\text{Tr}(f)\le 1$, and choose a projection $p\in N$ such that $\tau(p) = \text{Tr}(f)$. Then $\text{Tr}(pe_N) = \tau(p)$, so $f$ and $pe_N$ are equivalent projections in the factor $\langle M,e_N\rangle$. Let $V\in\langle M,e_N\rangle$ be a partial isometry with $VV^* = pe_N$ and $V^*V = f$. For $x,y\in N$,
\begin{equation}\label{SWWeq3.15}
VxV^* VyV^* = VxfyV^* = VfxyV^* = VxyV^*,
\end{equation}
and this shows that $x\mapsto VxV^*$ defines a $*$-homomorphism $\psi\colon \ N\to \langle M,e_N\rangle$. Since $e_NV=V$, the range of $\psi$ is contained in $e_N\langle M,e_N\rangle e_N = Ne_N$, and so there is a $*$-homomorphism $\phi\colon \ N\to N$ such that
\begin{equation}\label{SWWeq3.16}
VxV^* = \phi(x)e_N=e_N\phi(x),\qquad x\in N.
\end{equation}
If we multiply \eqref{SWWeq3.16} on the left by $V^*$ and use $f\in N'\cap \langle M,e_N\rangle$, $fV^*=V^*$ and $V^*e_N = V^*$, we obtain
\begin{equation}\label{SWWeq3.17}
xV^* = V^*\phi(x),\qquad x\in N.
\end{equation}
As in \cite{Popa.Betti} (see also the discussion in \cite[Section 8.4]{Sinclair.MasaBook}) there exists a non--zero partial isometry $w\in M$ such that
\begin{equation}\label{SWWeq3.18}
xw^* = w^*\phi(x),\qquad x\in N.
\end{equation}
Multiplication by $w$ on the right in \eqref{SWWeq3.18} shows that $w^*w\in N'\cap M = {\bb C}1$, and so $w$ is a unitary. Since $\phi(1)e_N = VV^* = pe_N$ from \eqref{SWWeq3.16}, while $\phi(1)=1$ from \eqref{SWWeq3.18}, it follows that $p=1$ and $\text{Tr}(f)=1$. Thus $\text{Tr}(f)<1$ is impossible, establishing (ii).

From \eqref{SWWeq3.18}, $wNw^* = \phi(N)\subseteq N$, and so $w\in \cl{ON}(N)$. Now consider $W = e_Nw\in \langle M,e_N\rangle$. This is a partial isometry because $WW^* = e_N$. Moreover, $W^*W  = w^*e_Nw\in N'\cap \langle M,e_N\rangle$ and is central by Lemma \ref{SWWlem3.2}. For $x\in N$,
\begin{align}
W^*Vx &= W^*\phi(x)V = w^*e_N\phi(x)V\notag\\
&= w^*\phi(x) e_NV = \label{SWWeq3.19} xw^*e_NV = xW^*V,
\end{align}
using \eqref{SWWeq3.17} and \eqref{SWWeq3.18}. Thus the operator $W^*V$ lies in $N'\cap \langle M,e_N\rangle$. Now
\begin{equation}\label{SWWeq3.20}
(W^*V)(W^*V)^* = W^*VV^*W = w^*e_Nw,
\end{equation}
so $W^*V$ is a partial isometry. Also
\begin{equation}\label{SWWeq3.21}
(W^*V)^*(W^*V) = V^*WW^*V = fV^*WW^*Vf
\end{equation}
and so is a projection $q$ in $N'\cap \langle M,e_N\rangle$ below $f$. From \eqref{SWWeq3.20}, $q$ is equivalent in $N'\cap\langle M,e_N\rangle$ to the central projection $w^*e_Nw$, so equality must hold. Thus $w^*e_Nw = q\le f$, and faithfulness of the trace Tr gives $w^*e_Nw = f$ since both projections have unit trace. This completes the proof of (i), and (ii) has already been proved.
\end{proof}

We will use this theorem in one of the forms below. Before we state these, we recall some simple facts from subfactor theory.

If $N\subseteq M$ is an irreducible inclusion of $\text{II}_1$ factors and $P$ is an intermediate von Neumann subalgebra, then any central element of $P$ lies in $N'\cap M = {\bb C}1$, showing that $P$ is automatically a subfactor of $M$. Moreover, the operator $e_P$ lies in $\langle M,e_N\rangle$. To see this, it suffices to show that $e_P$ commutes with $\langle M,e_N\rangle' = JNJ$. For each $x\in M$ and $y\in N$
\begin{equation}\label{SWWeq3.22}
(e_PJyJ - JyJe_P) x\xi = ({\bb E}_P(xy^*) - {\bb E}_P(x)y^*)\xi = 0
\end{equation}
since $y^*\in N\subseteq P$ so that ${\bb E}_P(xy^*) = {\bb E}_P(x)y^*$. A similar calculation also shows that $e_P\in N'$, so $e_P\in N'\cap \langle M,e_N\rangle$.

\begin{cor}\label{SWWcor3.5}
Let $N\subseteq M$ be an irreducible inclusion of $\mathrm{II}_1$ factors and let $P$ be an intermediate subfactor. Then $\cl{ON}(N)\subseteq P$ if and only if every nonzero projection $f\in N'\cap \langle M,e_N\rangle$ such that $f\le 1-e_P$ satisfies ${\rm Tr}(f)>1$.
\end{cor}

\begin{proof}
Suppose first that $\cl{ON}(N)\subseteq P$. Consider a nonzero projection $f\in N'\cap \langle M,e_N\rangle$ with $f\le 1-e_P$ and $\text{Tr}(f)\le 1$, a situation which will be shown to be impossible. By Theorem \ref{SWWthm3.4}, there exists $u\in \cl{ON}(N)$ such that $f=u^*e_Nu$. Since $u,u^*\in P$ and $e_Ne_P = e_N$, we obtain
\begin{equation}\label{SWWeq3.23}
e_Pfe_P = e_Pu^*e_Pe_Ne_Pue_P = u^*e_Pe_Ne_P u = u^*e_Nu = f.
\end{equation}
But $e_Pfe_P = 0$, so \eqref{SWWeq3.23} gives the contradiction $f=0$. Thus each such nonzero projection has $\text{Tr}(f)>1$.

Now suppose that every nonzero projection $f\in N'\cap \langle M,e_N\rangle$ with $f\le 1-e_P$ has $\text{Tr}(f)>1$. Consider an arbitrary unitary $u\in \cl{ON}(N)$. By Theorem \ref{SWWthm3.4}, $u^*e_Nu$ is a minimal projection in $N'\cap \langle M,e_N\rangle$ and is also central. Thus one of the inequalities $u^*e_Nu\le e_P$ or $u^*e_Nu \le 1-e_P$ must hold. However $\text{Tr}(u^*e_Nu) = 1$, so the second possibility is ruled out by the hypothesis. Thus $u^*e_Nu\le e_P$, from which it follows that $(1-e_P)u^*e_N = 0$. Apply this equation to the vector $\xi$ to obtain $u^*\xi = e_Pu^*\xi = {\bb E}_P(u^*)\xi$. Thus $u^*\in P$ so $u\in P$, and it follows that $\cl{ON}(N)\subseteq P$.
\end{proof}

\begin{cor}\label{SWWcor3.6}
Let $N\subseteq M$ be a finite index irreducible inclusion of ${\rm II}_1$ factors and let $P$ be an intermediate subfactor. Then $\cl{N}(N)\subseteq P$ if and only if every nonzero projection $f\in N'\cap \langle M,e_N\rangle$ with $f\le 1-e_P$ satisfies ${\rm Tr}(f)>1$.
\end{cor}

\begin{proof}
Under the additional hypothesis of finite index, $\cl{N}(N) = \cl{ON}(N)$ and the result follows from Corollary \ref{SWWcor3.5}.
\end{proof}

\begin{cor}\label{SWWcor3.7}
Let $N\subseteq M$ be a finite index irreducible inclusion of subfactors. Then $N$ is regular if and only if every minimal projection $f\in N'\cap \langle M,e_N\rangle$ satisfies ${\rm Tr}(f)=1$.
\end{cor}

\begin{proof}
Suppose that $N$ is not regular, and let $P = \cl{N}(N)''$, a proper intermediate subfactor. Then $1-e_P\ne 0$, so dominates a nonzero minimal projection $f\in N'\cap \langle M,e_N\rangle$ since this algebra is finite dimensional. By Corollary \ref{SWWcor3.6}, $\text{Tr}(f)>1$.

Conversely, suppose that $N$ is regular, so that $P$, as defined above, equals $M$. By Theorem \ref{SWWthm3.4} (ii), each nonzero projection $f\in N'\cap \langle M,e_N\rangle$ satisfies $\text{Tr}(f)\ge 1$. Suppose that there is a minimal $f$ with $\text{Tr}(f)>1$. Each $u\in \cl{N}(N)$ gives rise to a minimal projection $u^*e_Nu\in N'\cap \langle M,e_N\rangle$ and it is central, by Lemma \ref{SWWlem3.2}. Thus $fu^*e_Nu= 0$, otherwise this projection would contradict the minimality of $f$. If this equation is applied to the vector $u^*\xi$, then the result is
\begin{equation}\label{SWWeq3.24}
fu^*\xi = 0,\qquad u\in \cl{N}(N).
\end{equation}
Since $M$ is the strong closure of $\text{span}\{u^*\colon \ u\in \cl{N}(N)\}$, \eqref{SWWeq3.24} implies that $f$ annihilates $M\xi$ and so $f=0$, a contradiction which concludes the proof.
\end{proof}

\begin{rem}\label{SWWrem3.8}
We can apply the results above to a regular irreducible subfactor and recover the results of Jones, \cite{Jones.RegularSubfactors}. For a regular irreducible inclusion $N\subseteq M$ of $\text{II}_1$ factors, let $u_1,u_2,\dots\in\cl{N}(N)$ be left-coset representatives of $\cl{N}(N)/\cl{U}(N)$, which is a countable discrete group (see \cite[Section 3]{Jones.PropertiesMasas}), so that $(u_n^*e_Nu_n)_{n\geq 1}$ is a family of orthogonal projections in $N'\cap \langle M,e_N\rangle$, by Lemmas \ref{SWWlem3.2} and \ref{SWWlem3.3}. As in the proof of Corollary \ref{SWWcor3.7}, the projection $(\sum_{n\geq 1}u_n^*e_Nu_n)^\perp$ annihilates $M\xi$ so that 
$$
\sum_{n\geq 1}u_n^*e_Nu_n=1.
$$
Hence $(u_n)_{n\geq 1}$ forms a Pimsner-Popa basis for $N\subseteq M$, the index $[M:N]=\text{Tr}(1)$ is either an integer or $\infty$, and any one-sided normalizer of $N$ is automatically a two-sided normalizer of $N$.
\end{rem}

\section{Tensor products}\label{SWWsec4}

\indent 
Throughout this section, we will consider two irreducible inclusions $N_i \subseteq M_i$, $i=1,2$, of $\text{II}_1$ factors. Our objective is to relate the normalizer of the tensor product to the normalizers of the individual algebras. In the context of masas $A_i\subseteq M_i$, $i=1,2$, Chifan, \cite{Chifan.Normalisers}, has shown that $\cl{N}(A_1\ \ovl\otimes\ A_2)'' = \cl{N}(A_1)''\ \ovl\otimes\ \cl{N}(A_2)''$, and we will obtain a similar relationship for the $N_i$'s below. We will also be able to identify explicitly the normalizing unitaries for the tensor product.

We let $M = M_1\ \ovl\otimes\ M_2$ and $N=N_1\ \ovl\otimes\ N_2$. Tomita's commutant theorem ensures that $N \subseteq M$ is an irreducible inclusion. The basic construction behaves well with respect to tensor products, and there is a natural isomorphism
\begin{equation}\label{SWWeq4.1}
\langle M,e_N\rangle \cong \langle M_1,e_{N_1}\rangle\ \ovl\otimes\ \langle M_2,e_{N_2}\rangle
\end{equation}
where $e_N = e_{N_1} \otimes e_{N_2}$ (see Lemma \ref{SWWlem2.1}). If $\text{Tr}_i$, $i=1,2$, denote the canonical traces on $\langle M_i,e_{N_i}\rangle$, then the canonical trace Tr on $\langle M,e_N\rangle$ is given by $\text{Tr}= \text{Tr}_1 \otimes \text{Tr}_2$.

We begin by using the results of Section \ref{SWWsec3} to determine the one-sided normalizers of the tensor product.

\begin{thm}\label{SWWthm4.1}
Each unitary $v\in \cl{ON}(N)$ has the form $w(u_1\otimes u_2)$ where $w\in {\cl U}(N_1\ \ovl\otimes\ N_2)$ and $u_i\in \cl{ON}(N_i)$, $i=1,2$.
\end{thm}

\begin{proof}
It is clear that any unitary of the stated form is a one-sided normalizer of $N_1\ \ovl\otimes\ N_2$. Conversely, consider a one-sided unitary normalizer $v$ of $N_1\ \ovl\otimes\ N_2$. Then $v^*e_Nv$ is both a minimal and central projection in $N'\cap \langle M,e_N\rangle$, by Lemma \ref{SWWlem3.2}. Two applications of Tomita's commutant theorem show that
\begin{equation}\label{SWWeq4.2}
N'\cap \langle M,e_N\rangle = (N'_1 \cap \langle M_1,e_{N_1}\rangle)\ \ovl\otimes\ (N'_2\cap \langle M_2,e_{N_2}\rangle)
\end{equation}
and
\begin{equation}\label{SWWeq4.3}
Z(N'\cap\langle M,e_N\rangle) = Z(N'_1 \cap \langle M_1,e_{N_1}\rangle)\ \ovl\otimes\ Z(N'_2 \cap\langle M_2,e_{N_2}\rangle),
\end{equation}
where $Z(\cdot)$ denotes the center of an algebra. If these  centers are decomposed as direct sums of their atomic and diffuse parts, then minimality for $v^*e_Nv$ implies that $v^*e_Nv = p_1\otimes p_2$ for minimal projections $p_i\in Z(N'_i\cap \langle M_i,e_{N_i}\rangle)$, $i=1,2$. By Theorem \ref{SWWthm3.4}, $\text{Tr}_i(p_i)\ge 1$, forcing equality since $\text{Tr}(v^*e_Nv) = 1$. A second application of Theorem \ref{SWWthm3.4} gives the existence of unitaries $u_i\in \cl{ON}(N_i)$ such that $p_i = u^*_ie_{N_i}u_i$ for $i=1,2$. Thus
\begin{equation}\label{SWWeq4.4}
v^*e_Nv = u^*_1e_{N_1}u_1 \otimes u^*_2e_{N_2}u_2 = (u_1\otimes u_2)^* e_N(u_1\otimes u_2).
\end{equation}
By Lemma \ref{SWWlem3.3}, there exists a unitary $w\in {\cl U}(N_1\ \ovl\otimes\ N_2)$ such that $v=w(u_1 \otimes u_2)$.
\end{proof}

The case of one-sided normalizers above easily leads to a similar result for unitary normalizers.

\begin{thm}\label{SWWthm4.2}
Each unitary $v\in \cl{N}(N)$ has the form $w(u_1\otimes u_2)$ where $w\in {\cl U}(N_1\ \ovltimes\ N_2)$ and $u_i\in \cl{N}(N_i)$, $i=1,2$.
\end{thm}

\begin{proof}
Clearly each unitary of the stated form is a unitary normalizer of $N$. Conversely, let $v\in \cl{N}(N)$. Viewing $v$ as a one-sided normalizer, Theorem \ref{SWWthm4.1} implies that $v$ has the form $v=w(u_1\otimes u_2)$ where $w\in {\cl U}(N_1\ \ovltimes\ N_2)$ and $u_i \in \cl{ON}(N_i)$ for $i=1,2$. Then $w^*v\in \cl{N}(N_1\ \ovltimes\ N_2)$, showing that $u_i\in \cl{N}(N_i)$, otherwise $(u_1\otimes u_2)^*(N_1\ \ovltimes\ N_2)(u_1\otimes u_2)$ would be strictly contained in $N_1\ \ovltimes\ N_2$.
\end{proof}

\begin{cor}\label{SWWcor4.3}
Let $N_i\subseteq M_i$, $i=1,2$, be irreducible inclusions of ${\rm II}_1$ factors. Then
\begin{equation}\label{SWWeq4.5}
\cl{ON}(N_1\ \ovltimes\ N_2)'' = \cl{ON}(N_1)'' ~\ovltimes~ \cl{ON}(N_2)''
\end{equation}
and
\begin{equation}\label{SWWeq4.6}
\cl{N}(N_1\ \ovltimes\ N_2)'' = \cl{N}(N_1)'' ~\ovltimes~ \cl{N}(N_2)''.
\end{equation}
\end{cor}

\begin{proof}
This is immediate from the characterizations of one-sided and two-sided normalizers in Theorems \ref{SWWthm4.1} and \ref{SWWthm4.2}.
\end{proof}

\begin{cor}\label{SWW.cor4.4}
Let $N_i\subseteq M_i$, $i=1,2$ be inclusions of ${\rm II}_1$ factors with $N_i$ singular in $M_i$. Then $N_1\ \ovl\otimes\ N_2$ is singular in $M_1\ \ovl\otimes\ M_2$.
\end{cor}
\begin{proof}
This is immediate from Corollary \ref{SWWcor4.3} as a singular subfactor is automatically irreducible.
\end{proof}

\begin{rem}\label{SWWrem4.5}
The result of Theorem \ref{SWWthm4.2} can be false in more general situations. Consider two regular masas $A_i\subseteq M_i$, $i=1,2$. Then $A_1\ \ovltimes\  A_2$ is a regular masa in $M_1\ \ovltimes\ M_2$ since
\[
\{u_1\otimes u_2\colon \ u_i\in {\cl N}(A_i)\} \subseteq {\cl N}(A_1\ \ovltimes\ A_2).
\]
Choose projections $p_i\in A_i$ of equal trace 1/2, so that $p_1\otimes 1$ and $1\otimes p_2$ have equal trace in $A_1\ \ovltimes\ A_2$. From \cite{Popa.NotesCartan}, there exists a unitary $u\in {\cl N}(A_1\ \ovltimes\ A_2)$ such that $u(p_1\otimes 1)u^* = 1\otimes p_2$. Then $u$ cannot have the form $(u_1\otimes u_2)w$ for $w\in {\cl U}(A_1\ \ovltimes\ A_2)$ and $u_i\in {\cl N}(A_i)$, since this would imply that
\begin{align}
1\otimes p_2 &= u(p_1\otimes 1)u^* = (u_1\otimes u_2) (p_1\otimes 1)(u_1\otimes u_2)^*\notag\\
\label{SWWeq4.7}
&=  (u_1p_1u^*_1\otimes 1),
\end{align}
which is impossible.
\end{rem}

\section{Fixed point algebras}\label{SWWsec5}

\indent 

In this section we consider a finite group $G$ acting by outer automorphisms on a $\text{II}_1$ factor $M$, where we denote the action by the homomorphism $\alpha\colon \ G\to \text{Aut}(M)$. For each subgroup $H\subseteq G$, the fixed point algebra is denoted by
\[
M^H = \{x\in M\colon \ \alpha_h(x) = x \quad \text{for}\quad h\in H\}.
\]
The following lemma is quoted from \cite{Nakamura.GaloisFinite} (see also \cite{Jones.SubfactorsBook}) for the reader's convenience.

\begin{lem}\label{SWWlem5.1}
Let $G$ be a finite group, let $M$ be a ${\rm II}_1$ factor, and let $\alpha\colon \ G\to {\rm Aut}(M)$ be an action of $G$ on $M$ by outer automorphisms. Then $M^G$ is a factor and the inclusion $M^G\subseteq M$ is irreducible.
\end{lem}

Let $N=M^G$. For each $g\in G$, define a unitary $u_g\in B(L^2(M))$ by
\begin{equation}\label{SWWeq5.6}
u_g(x\xi) = \alpha_g(x)\xi,\qquad x\in M.
\end{equation}
Then $u_gxu^*_g = \alpha_g(x)$ for $x\in M$, and the crossed product $M\rtimes_\alpha G$ is faithfully represented on $L^2(M)$ as the set of operators of the form $\sum_{g\in G}x_gu_g$ for $x_g\in M$ and $g\in G$. A simple calculation identifies the commutant as $JM^GJ = JNJ$, so that in this representation $M\rtimes_\alpha G$ is identified with $\langle M,e_N\rangle$. The projection $e_N$ is $\frac1{|G|} \sum\limits_{g\in G}u_g$. From \cite{Jones.Index},
\begin{equation}\label{SWWeq5.7}
[M:N] = [\langle M,e_N\rangle:M] = [M\rtimes_\alpha G:M] = |G|,
\end{equation}
and so the canonical trace Tr on $\langle M,e_N\rangle$ satisfies $\text{Tr}(1)=|G|$. Its restriction to $M\rtimes_\alpha G$ is the unique trace on this factor (up to scaling), and so $\text{Tr}(u_g) = 0$ for $g\ne e$. Consequently Tr agrees with the trace arising from $x\mapsto |G| \langle x\delta_e,\delta_e\rangle$ on the group algebra ${\bb C}G$. Moreover, $N'\cap \langle M,e_N\rangle$ is $\text{span}\{u_g\colon \ g\in G\} = {\bb C}G$ since any element $\sum\limits_{g\in G} x_gu_g\in N'\cap \langle M,e_N\rangle$ satisfies
\begin{equation}\label{SWWeq5.8}
\sum_{g\in G} vx_gv^*u_g = \sum_{g\in G}x_gu_g,\qquad v\in {\cl U}(N),
\end{equation}
implying that each $x_g$ is a scalar by Lemma \ref{SWWlem5.1}.

\begin{lem}\label{SWWlem5.2}
Suppose that $G$ is a finite abelian group acting on a ${\rm II}_1$ factor $M$ by outer automorphisms $\alpha\colon\ G\to {\rm Aut}(M)$. Then $M^G$ is regular in $M$.
\end{lem}

\begin{proof}

For a general finite group $G$, the minimal central projections in the group algebra $\mathbb CG$ correspond to irreducible representations of $G$. Furthermore the trace of the minimal central projection corresponding to a $d$-dimensional representation is $d^2$. By Schur's Lemma, every irreducible representation of an abelian group is one-dimensional. As $\mathbb CG$ is also abelian, the minimal projections in $\mathbb CG$ are central and so have trace $1$. The result then follows from Corollary \ref{SWWcor3.7}.
\end{proof}

\begin{rem}
In the situation of Lemma \ref{SWWlem5.2}, the minimal projections in $\mathbb CG$ are of the form $e_\phi = \frac1{|G|} \sum\limits_{g\in G} \phi(g)u_g$ for $\phi$ in the dual group $\hat{G}$.
\end{rem}

Before we can state the main result of this section, Theorem \ref{SWWthm5.4}, we require a technical result on actions of quotient groups. If $G$ is a finite group with a normal subgroup $H$, then any action $\alpha\colon \ G\to\text{Aut}(M)$ induces an action $\theta\colon \ G/H\to \text{Aut}(M^H)$ by
\begin{equation}\label{SWWeq5.9}
\theta_{gH}(x) = \alpha_g(x), \qquad x\in M^H,\quad g\in G.
\end{equation}

\begin{lem}\label{SWWlem5.3}
Let $G$ be a finite group acting on a ${\rm II}_1$ factor $M$ by outer automorphisms $\alpha\colon \ G\to {\rm Aut}(M)$ and let $H$ be a normal subgroup. Then the action $\theta\colon \ G/H\to {\rm Aut}(M^H)$ is also by outer automorphisms.
\end{lem}

\begin{proof}
Fix $g\in G\backslash H$, and fix an ordering $g_1,\ldots, g_n$ of the elements of $G$ so that $g_1,\ldots, g_p$ is a list of the elements of $H$. Left multiplication by $g$ induces a permutation $\pi$ of the integers $\{1,\ldots, n\}$ such that $gg_i = g_{\pi(i)}$, $1\le i\le n$. 

Suppose that $\theta_{gH}$ is implemented by a unitary $v\in M^H$. This will be shown to be impossible by deriving the contradiction $v=0$. Note that the elements of $M^H$ have the form $\sum\limits^p_{i=1} \alpha_{g_i}(x)$ for $x\in M$. Applying $\theta_{gH}$ gives
\begin{equation}\label{SWWeq5.10}
\theta_{gH} \left(\sum^p_{i=1} \alpha_{g_i}(x)\right) v = v \sum^p_{i=1} \alpha_{g_i}(x),\qquad x\in M,
\end{equation}
and so
\begin{equation}\label{SWWeq5.11}
\sum^p_{i=1} \alpha_{g_{\pi(i)}}(x) v = v \sum^p_{i=1} \alpha_{g_i}(x),\qquad x\in M.
\end{equation}
Multiplication by an arbitrary element $y\in M'$ in \eqref{SWWeq5.11} leads to
\begin{equation}\label{SWWeq5.12}
\sum^p_{i=1} y \alpha_{g_{\pi(i)}}(x)v = v \sum^p_{i=1} y \alpha_{g_i}(x),\qquad x\in M,\quad y\in M'.
\end{equation}
We now introduce some matrices as follows, using $T$ to denote the transpose:\ 
$$
R = (\underbrace{1,1,\ldots, 1}_{n}),\quad C=R^T,
$$
\begin{alignat*}{2}
V_1 &= (v_1,\ldots, v_n)^T\quad &\text{where}\quad v_j &= \begin{cases}
v,&\text{$j\in\{\pi(1),\ldots,\pi(p)\}$,}\\
0,&\text{otherwise,}
                                                          \end{cases}\\
V_2 &= (\tilde v_1,\ldots, \tilde v_n) \quad &\text{where}\quad \tilde v_j &= \begin{cases}
v,&\text{$1\le j\le p$,}\\
0,&\text{otherwise,}
                                                         \end{cases}
\end{alignat*}
and
\[
D(x,y) = \begin{pmatrix}
y\alpha_{g_1}(x)\\ &\ddots\\ &&y\alpha_{g_n}(x)
          \end{pmatrix}.
\]
 With this notation, \eqref{SWWeq5.12} becomes
\begin{equation}\label{SWWeq5.13}
RD(x,y)V_1 = V_2 D(x,y)C,\qquad x\in M,\quad y\in M'.
\end{equation}
By Lemma \ref{SWWlem2.2}, the operators $D(x,y)$ generate the von Neumann algebra 
\begin{equation}\label{SWWeq5.14}
\left\{\begin{pmatrix} t_1\\ &\ddots\\ &&t_n\end{pmatrix}\colon \ t_i \in 
B(L^2(M))\right\}.
\end{equation}
Thus $D(x,y)$ can be replaced in \eqref{SWWeq5.13} by any diagonal matrix from 
\eqref{SWWeq5.14}, leading to 
\begin{equation}\label{SWWeq5.15}
\sum^p_{i=1} t_{\pi(i)}v = \sum^p_{i=1} vt_i,\qquad t_i\in B(L^2(M)).
\end{equation}
From this we deduce that $v$ is a scalar by  setting $t_1=t_2 = \cdots = t_n = t$,
where $t\in B(L^2(M))$ is arbitrary. If $v\ne 0$, then it can be canceled from \eqref{SWWeq5.15}. Consequently
\begin{equation}\label{SWWeq5.16}
\sum^p_{i=1} t_{\pi(i)} = \sum^p_{i=1} t_i,\qquad t_i\in B(L^2(M)).
\end{equation}
Since $g\notin H$, $\{{\pi(1)},\ldots,{\pi(p)}\}\cap\{1,\ldots, p\}=\emptyset$, and so \eqref{SWWeq5.16} can fail for suitable choices of the operators $t_j$, $1\le j \leq n$. This forces the contradiction $v=0$, showing that each $\theta_{gH}$ for $g\in G\backslash H$ is outer.
\end{proof}

We are now able to describe the algebra ${\cl N}(M^G)''$. We use the standard notation $[G,G]$ for the commutator subgroup
\[
[G,G] = \langle g^{-1}h^{-1}gh\colon \ g,h\in G\rangle.
\]
This is the smallest normal subgroup for which the quotient is abelian, and can also be described as
\begin{equation}\label{SWWeq5.17}
[G,G] = \cap\{H\triangleleft G\colon \ G/H \text{ is abelian}\}.
\end{equation}

\begin{thm}\label{SWWthm5.4}
Let $G$ be a finite group acting on a ${\rm II}_1$ factor $M$ by outer automorphisms $\alpha\colon \ G\to {\rm Aut}(M)$. Then
\[
{\cl N}(M^G)'' = M^{[G,G]}.
\]
\end{thm}

\begin{proof}
By the Galois theory of \cite{Nakamura.GaloisFinite,Nakamura.Galois}, there is a subgroup $K$ such that ${\cl N}(M^G)'' = M^K$, and so it must be shown that $K = [G,G]$.

First suppose that $H$ is a normal subgroup of $G$ with abelian quotient $G/H$. There is a well defined action of $G/H$ on $M^H$ given by
\begin{equation}\label{SWWeq5.18}
\theta_{gH}(x) = \alpha_g(x),\qquad x\in M^H,\quad g\in G. 
\end{equation}
By Lemma \ref{SWWlem5.3}, this is an outer action and $G/H$ is abelian. By Lemma \ref{SWWlem5.2}, $(M^H)^{G/H}$ is regular in $M^H$ and the first of these algebras is $M^G$. Consequently ${\cl N}(M^G)''$ contains $M^H$. From the Galois theory of \cite{Nakamura.GaloisFinite,Nakamura.Galois}, $K\subseteq H$. It then follows from \eqref{SWWeq5.17} that $K\subseteq [G,G]$.

Conversely, suppose that $u$ is a unitary normalizer of $M^G$. For each $x\in M^G$, $uxu^*\in M^G$ and so
\begin{equation}\label{SWWeq5.19}
uxu^* = \alpha_g(uxu^*) = \alpha_g(u)x\alpha_g(u)^*,\qquad g\in G.
\end{equation}
Then $u^*\alpha_g(u)$ lies in the relative commutant of $M^G$ which is trivial by Lemma \ref{SWWlem5.1}. Thus there exists a homomorphism $g\mapsto \lambda_g$ of $G$ into the circle group ${\bb T}$ such that $\alpha_g(u) = \lambda_gu$. The kernel $H$ of this homomorphism is a normal subgroup with abelian quotient and $u\in M^H$. Since $[G,G]\subseteq H$, it follows that $u\in M^{[G,G]}$. Thus $M^K \subseteq M^{[G,G]}$, and we obtain the reverse inclusion $[G,G] \subseteq K$, proving equality.
\end{proof}

\begin{exm}
We are able to obtain more examples of singular subfactors from Theorem \ref{SWWthm5.4}. If $G$ is any simple non-abelian finite group, such as $A_5$, acting by outer automorphisms on a ${\rm II}_1$ factor $M$, then $M^G$ is singular in $M$. We can also obtain singular subfactors from outer actions of non-simple groups.  Consider $G=H\times H$, where $H$ is a finite simple non-abelian group. The normal subgroups of $G$ are $\{e\}\times \{e\},\ H\times \{e\}$,  $\{e\}\times H$, and $H\times H$, so $[G,G]=G$ and $M^G$ is singular in $M$ whenever $G$ acts on the ${\rm II}_1$ factor $M$ by outer automorphisms.  It may be of interest to examine the strong-singularity constants of \cite{Wiggins.StrongSing} arising from these fixed point examples.
\end{exm}

\section{Group factors}\label{SWWsec6}

\indent

This section is concerned with inclusions of $\text{II}_1$ factors which arise from inclusions of countable discrete groups. We examine the normalizers of $L(H) \subseteq L(G)$ when $H\subseteq G$, and relate these to the algebraic normalizers of $H$ as a subgroup of $G$. One-sided normalizers will again play a role, so apart from the standard notation
\[
{\cl N}_G(H) = \{g\in G\colon \ gHg^{-1}=H\}
\]
for the normalizer, we also introduce the semigroup of one-sided normalizers
\[
\cl{ON}_G(H) = \{g\in G\colon \ gHg^{-1} \subseteq H\}.
\]
Subsequently we will exhibit situations where these two normalizers are distinct.

We will need the following simple lemma, where some of the conditions are reminiscent of those in \cite{MvN.4,Dixmier.Masa}, concerning masas arising from subgroups. We will use the terminology $K$-\emph{conjugates of a group element} $g\in G$ for $\{kgk^{-1}\colon \ k\in K\}$ where $K$ is a subgroup of $G$.

\begin{lem}\label{SWWlem6.1}
Let $K\subseteq H\subseteq G$ be an inclusion of countable discrete groups.
\begin{itemize}
\item[\rm (i)] $L(H)$ is irreducible in $L(G)$ if and only if each $g\in G\backslash\{e\}$ has infinitely many $H$-conjugates;
\item[\rm (ii)] If $G$ is I.C.C., $L(H)$ is irreducible in $L(G)$ and $K$ has finite index in $H$, then $L(K)$ is irreducible in $L(G)$.
\end{itemize}
\end{lem}

\begin{proof}
\n (i)~~If there exists $g\in G\backslash\{e\}$ with finitely many $H$-conjugates $\{h_igh^{-1}_i\colon \ 1\le i\le n\}$, then $\sum\limits^n_{i=1} \lambda_{h_i}\lambda_g \lambda_{h^{-1}_i}$ is a non-scalar element in $L(H)' \cap L(G)$.

Conversely, suppose that each $g\in G\backslash\{e\}$ has infinitely many $H$-conjugates, and let $y\in L(H)'\cap L(G)$. Viewing $y$ as an $\ell^2(G)$-vector $\sum\limits_{g\in G} \alpha_g\delta_g$ where $\alpha_g\in {\bb C}$ and $\sum\limits_{g\in G} |\alpha_g|^2 < \infty$, commutation with $h\in H$ implies that
\begin{equation}\label{SWWeq6.1}
\sum_{g\in G} \alpha_g\delta_{gh} = \sum_{g\in G} \alpha_g \delta_{hg} = \sum_{g\in G} \alpha_{h^{-1}gh} \delta_{gh}
\end{equation}
so that $\alpha_g = \alpha_{h^{-1}gh}$ for $g\in G$ and $h\in H$. Thus $\alpha_g=0$ for $g\ne e$ since otherwise there are infinitely many equal non-zero values for the coefficients. It follows that $y\in {\bb C}1$ and that $L(H)$ is irreducible in $L(G)$.

\n (ii) Write $H=\bigcup_{i=1}^mh_iK$ as a finite union of left $K$-cosets.  The $H$-conjugates of $g\in G$ with $g\neq e$ are
$$
\bigcup_{i=1}^m\{h_ikgk^{-1}h_i^{-1}\colon k\in K\}.
$$
By (i), the set of $H$-conjugates of $g$ is infinite. Hence the $K$-conjugates of $g$ are also infinite. Another application of (i) shows that $L(K)$ is irreducible in $L(G)$.
\end{proof}

\begin{thm}\label{SWWthm6.2}
Let $H\subseteq G$ be an inclusion of countable discrete groups, where $G$ is I.C.C.\ and $L(H)$ is irreducible in $L(G)$.
\begin{itemize}
\item[\rm (i)] Each $u\in \cl{ON}(L(H))$ has the form $w\lambda_g$ for $w\in {\cl U}(L(H))$ and $g\in \cl{ON}_G(H)$; 
\item[\rm (ii)] each $u\in \cl{N}(L(H))$ has the form $w\lambda_g$ for $w\in {\cl U}(L(H))$ and $g\in {\cl N}_G(H)$.
\end{itemize}
\end{thm}

\begin{proof}
\n (i)~~Let $u$ be a fixed element of $\cl{ON}(L(H))$, and write $v=u^*$ (just as in Section \ref{SWWsec3}, it is more convenient to consider $u^*$). Then $x\mapsto v^*xv$ defines a $*$-homomorphism $\alpha\colon \ L(H)\to L(H)$ satisfying $xv=v\alpha(x)$ for $x\in L(H)$. Let $\{g_iH\}^\infty_{i=1}$ be a listing of the left $H$-cosets in $G$ (which could also be a finite set if $H$ has finite index). The closed subspaces $\lambda_{g_i}\ell^2(H)$, $i\ge 1$, of $\ell^2(G)$ are pairwise orthogonal and span this Hilbert space. Viewing $v$ as an element of $\ell^2(G)$, we may then write $v = \sum\limits^\infty_{i=1} \lambda_{g_i}\eta_i$ for vectors $\eta_i\in \ell^2(H)$ satisfying $\sum\limits^\infty_{i=1} \|\eta_i\|^2_2 < \infty$. In fact, $\eta_i=\mathbb E_{L(H)}(\lambda_{g_i}^*v)\in L(G)$. This gives
\begin{equation}\label{SWWeq6.2}
\sum^\infty_{i=1} \lambda_{hg_i}\eta_i = \sum^\infty_{i=1} \lambda_{g_i}J\alpha(\lambda_h)^*J\eta_i,\qquad h\in H.
\end{equation}
By renumbering, we may assume that $\|\eta_1\|_2\ne 0$. Each $h\in H$ defines a permutation of the left $H$-cosets by $g_iH \mapsto hg_iH$, and so there is a permutation $\pi_h$ of $\{1,2,\ldots\}$ such that $hg_iH = g_{\pi_h(i)}H$. The map $h\mapsto \pi_h$ is then a homomorphism of $H$ into the group of permutations of ${\bb N}$. Moreover, there are maps $\phi_i\colon \ H\to H$ such that $hg_i = g_{\pi_h(i)}\phi_i(h)$, $h\in H$, and \eqref{SWWeq6.2} becomes
\begin{equation}\label{SWWeq6.2a}
\sum^\infty_{i=1} \lambda_{g_{\pi_{h}(i)}}\lambda_{\phi_i(h)}\eta_i = \sum^\infty_{i=1} \lambda_{g_i}J\alpha(\lambda_h)^*J\eta_i,\qquad h\in H.
\end{equation}

Define an equivalence relation on ${\bb N}$ by $i\sim j$ if $\|\eta_i\|_2 = \|\eta_j\|_2$, and let $S_i$ be the equivalence class containing $i$. For each $i\ge 1$, the $\ell^2(H)$-components of the vectors in \eqref{SWWeq6.2a} which lie in $\lambda_{g_{\pi_{h}(i)}}\ell^2(H)$ are $\lambda_{\phi_i(h)}\eta_i$ on the left and $J\alpha(\lambda_h)^*J \eta_{\pi_h(i)}$ on the right. Their norms are respectively $\|\eta_i\|_2$ and $\|\eta_{\pi_h(i)}\|_2$, showing that $i\sim \pi_h(i)$ for all $i\ge 1$ and $h\in H$. Since $\|\eta_1\|\ne 0$ and $\lim\limits_{i\to\infty} \|\eta_i\|_2= 0$, the set $S_1$ is finite and each $\pi_h$ maps $S_1$ to $S_1$. The restriction of $\pi_h$ to $S_1$ gives a homomorphism of $H$ into the finite group of permutations  of $S_1$, and so the kernel $K$ has finite index in $H$. Then $kg_1H = g_1H$ for $k\in K$, so there is a homomorphism $\beta\colon K \to H$ given by $\beta (k)=g_1^{-1}kg_1$. We also regard $\beta$ as being a $*$-automorphism of $L(G)$ defined by $\beta (x)=\lambda^*_{g_1}x\lambda_{g_1}$ and  which maps $L(K)$ to $L(H)$.
In \eqref{SWWeq6.2} replace $h\in H$ by $k\in K \subseteq
H$, multiply on the left by $\lambda^*_{g_1}$, and
compare the vectors that lie in $\ell^2(H)$. The result is
\begin{equation}\label{SWWeq6.3}
\beta(\lambda_k)\eta_1 = J{\alpha(\lambda^*_k)}J\eta_1,\qquad k\in K.
\end{equation}
Now $\eta_1 = {\bb E}_{L(H)}(\lambda^{*}_{g_1}v)$, and so $\eta_1$ can be viewed as an element $y\in L(H)$. Thus \eqref{SWWeq6.3} becomes
\begin{equation}\label{SWWeq6.4}
\beta(\lambda_k)y = y\alpha(\lambda_k),\qquad k\in K,
\end{equation}
and so $yy^* \in \beta(L(K))'\cap L(H)$. Thus $\lambda_{g_1}yy^*\lambda^*_{g_1}\in L(K)'\cap L(G)$.

By Lemma \ref{SWWlem6.1} (ii), $L(K)$ is irreducible in $L(G)$. Thus $\lambda_{g_1}yy^*\lambda^*_{g_1} \in {\bb C}1$, so by scaling we may assume that $y$ is a unitary $w^*\in L(H)$. Then \eqref{SWWeq6.4} becomes
\begin{equation}\label{SWWeq6.5}
v^*\lambda_kv = w\lambda^*_{g_1}\lambda_k\lambda_{g_1}w^*,\qquad k\in K,
\end{equation}
showing that $\lambda_{g_1}w^*v^*\in L(K)' \cap L(G) = {\bb C}1$. Multiplying $w^*$ by a suitable $e^{i\theta}$, we may assume that $\lambda_{g_1}w^*v^*=1$, so $u=v^*=w\lambda_{g^{-1}_1}$. Since $uL(H) u^* \subseteq L(H)$, we see that $\lambda_{g^{-1}_1}L(H) \lambda_{g_1} \subseteq L(H)$, so if we replace $g^{-1}_1$ by $g$ then $g\in \cl{ON}_G(H)$ and $u$ has the form $w\lambda_g$ for $w\in {\cl U}(L(H))$.

\n (ii)~~If $u\in {\cl N}(L(H))$ then $u$ has the form $w\lambda_g$ for $g\in \cl{ON}_G(H)$ and $w\in {\cl U}(L(H))$. Since conjugation by $u$ maps $L(H)$ onto $L(H)$, the same is true for conjugation by $\lambda_g$, showing that $gHg^{-1} = H$ and that $g\in {\cl N}_G(H)$.
\end{proof}

\begin{exm}\label{SWWexm6.3}
An immediate consequence of Theorem \ref{SWWthm6.2} (ii) is that $L(H)$ is singular in $L(G)$ precisely when ${\cl N}_G(H)=H$. Here we give examples of singular inclusions $L(H)\subseteq L(G)$ which nevertheless have non-trivial one-sided normalizers.

Consider the free group ${\bb F}_\infty$, where the generators are written $\{g_i\colon \ i\in {\bb Z}\}$, and for each $n\in {\bb Z}$, let $H_n$ be the subgroup generated by $\{g_i\colon \ i\ge n\}$. The shift $i\mapsto i+1$ on ${\bb Z}$ induces an automorphism $\phi$ of ${\bb F}_\infty$ defined on generators by $\phi(g_i) = g_{i+1}$, $i\in {\bb Z}$, and $\phi$ maps $H_n$ into $H_{n+1} \subseteq H_n$. Then $n\mapsto \phi^n$ gives a homomorphism $\alpha\colon \ {\bb Z}\to {\rm Aut}({\bb F}_\infty)$, and we let $G$ be the semidirect product ${\bb F}_\infty \rtimes_\alpha {\bb Z}$. We abuse notation and write the elements of this group as $w\phi^n$ where $w\in {\bb F}_\infty$. The multiplication is
\begin{equation}\label{SWWeq6.6}
(v\phi^n)(w\phi^m) = (v\phi^n(w))\phi^{n+m},\qquad v,w\in {\bb F}_\infty,\quad m,n\in {\bb Z}.
\end{equation}
We then consider the inclusion $H_n\subseteq {\bb F}_\infty \rtimes_\alpha {\bb Z}$. By construction $\phi H_n\phi^{-1} = H_{n+1}$ so $\phi$ is a one-sided normalizer of $H_n$ for each $n\in {\bb Z}$. We now show that the only normalizers of $H_n$ lie in $H_n$.

Suppose that $v\phi^k$ has the property that $v\phi^kH_n\phi^{-k}v^{-1} = H_n$. If $v\in H_n$ then $H_{n+k} = H_n$, forcing $k=0$, and we see that $v\phi^k\in H_n$. Thus we may assume that $v\notin H_n$. Let $j$ be the minimal integer such that $g_j$ appears in $v$. Then $j<n$ otherwise $v\in H_n$. Then $H_{n+k} = v^{-1} H_nv\subseteq H_j$, so $n+k\ge j$. Take $r>n$ such that the letter $g_r$ does not appear in the reduced word $v$. Then there is no cancellation in $v^{-1}g_rv$. In particular, the letter $g_j$ cannot cancel from $v^{-1}g_rv\in v^{-1}H_nv = H_{n+k}$, and so $n+k\le j$, showing that $v^{-1}H_nv = H_j$. There is also no cancellation in $vg_rv^{-1}$, so $vg_rv^{-1}\in vH_jv^{-1}$ is not contained in $H_n$, a contradiction. Thus there are no non-trivial normalizers of $H_n$, so $L(H_n)$ is singular in $L(G)$ although it does have non-trivial one-sided normalizers. Further algebraic calculations along the same lines show that
\begin{equation}\label{SWWeq6.7}
\cl{ON}_G(H_n) = \{v\phi^r\colon \ v\in H_n,\quad r\ge 0\},\qquad n\in {\bb Z}.
\end{equation}
We omit the easy details.

It is worth noting that the disparity between normalizers and one-sided normalizers
in this example is extreme; the former generate $L(H_n)$ while the latter generate
$L(G)$.
\end{exm}

\begin{rem}
Just as in Remark \ref{SWWrem4.5}, the analgous statement to Theorem \ref{SWWthm6.2} is false in the abelian situation. Let $H$ be an abelian subgroup of an I.C.C. group $G$ such that every element of $G\setminus H$ has infinitely many $H$-conjugates --- this is Dixmier's condition, \cite{Dixmier.Masa}, which is equivalent to $L(H)$ being a masa in $L(G)$. Normalizers of $L(H)$ are not necessarily of the form $u\lambda_g$ for some $g\in\cl{N}_G(H)$ and a unitary $u\in L(H)$. This leads to a question to which we do not know the answer. Suppose that $\cl{N}_G(H)=H$.  Must $L(H)$ be singular in $L(G)$?  The methods used to prove singularity of masas coming from subgroups, \cite{Sinclair.StrongSing}, \cite[Lemma 2.1]{Sinclair.StrongSing2}, require additional algebraic conditions on $H\subseteq G$.
\end{rem}

\section{Bimodules}\label{SWWsec7}

\indent 

We consider an inclusion $H\subseteq G$ of countable discrete groups such that $L(H)' \cap L(G) = {\bb C}1$. Each projection $f\in L(H)'\cap \langle L(G), e_{L(H)}\rangle$ has a range which is invariant under left and right multiplications by elements of $L(H)$. Conversely, any norm closed $L(H)$-bimodule in $\ell^2(G)$ is the range of such a projection $f$. In this section we investigate the structure of $L(H)'\cap \langle L(G), e_{L(H)}\rangle$. In Theorem \ref{SWWthm7.1} we characterize the bimodules for projections of finite trace, and subsequently we show that the structure can be much more complicated when projections of infinite trace are considered. We recall from Theorem \ref{SWWthm3.4} (ii) that any non--zero projection $f\in L(H)'\cap \langle L(G),e_{L(H)}\rangle$ satisfies $\text{Tr}(f)\ge 1$.

\begin{thm}\label{SWWthm7.1}
Let $H\subseteq G$ be an inclusion of countable discrete groups such that $L(H)' \cap L(G) = {\bb C}1$, and let $f\in L(H)'\cap \langle L(G), e_{L(H)}\rangle$ be a non--zero projection such that $\mathrm{Tr}(f)<\infty$. Then $\mathrm{Tr}(f)$ is an integer, and there exist $g_1,\ldots, g_n\in G$ such that the range of $f$ is the direct sum $\bigoplus\limits^n_{i=1} \lambda_{g_i}\ell^2(H)$. In particular, the range of $f$ is a finite sum of $L(H)$-bimodules generated by double cosets $HgH$ each of which is a finite sum of right $L(H)$-modules generated by left cosets $gH$.
\end{thm}

\begin{proof}
The projection onto $\bigoplus\limits^n_{i=1} \lambda_{g_i}\ell^2(H)$ is $\sum\limits^n_{i=1} \lambda_{g_i}e_{L(H)}\lambda^{*}_{g_i}$, and so has trace equal to $n$. We consider a  non--zero projection $f\in L(H)'\cap \langle L(G), e_{L(H)}\rangle$ with $\text{Tr}(f)<\infty$, and we write $\text{Tr}(f) = (n-1)+\mu$ where $n\in {\bb N}$ and $\mu\in (0,1]$. In the course of the proof it will be shown that $\mu=1$.

Choose a projection $p\in L(H)$ with $\tau(p)=\mu$. Following the approach of \cite{Popa.Entropy}, the diagonal projections
\begin{equation}\label{SWWeq7.1}
P_1 = \begin{pmatrix}
f\\ &0\\ &&\ddots\\ &&&0\\ &&&&0
      \end{pmatrix}\quad \text{and}\quad P_2 = 
\begin{pmatrix}
e_{L(H)}\\ &e_{L(H)}\\ &&\ddots\\ &&&e_{L(H)}\\ &&&&pe_{L(H)}
\end{pmatrix}
\end{equation}
in ${\bb M}_n(\langle L(G), e_{L(H)}\rangle)$ have equal finite traces and so are equivalent in this factor. Thus there exists a column matrix $V = (v_1,\ldots, v_n)^T$ with entries $v_i\in \langle L(G),e_{L(H)}\rangle$ such that $V^*V=f$ and $VV^* = P_2$. In particular, $v^*_ie_{L(H)} = v^*_i$, $1\le i\le n$. As in \cite{Popa.Betti} (see also \cite[Section 8.4]{Sinclair.MasaBook}), the map $x\mapsto VxV^*$ defines a homomorphism $\psi\colon \ L(H)\to {\bb M}_n\langle L(G),e_{L(H)}\rangle$ whose range lies under $P_2$, and so there is a homomorphism $\phi\colon \ L(H)\to {\bb M}_n(L(H))$ such that 
\begin{equation}\label{SWWeq7.2}
\psi(x) = \phi(x)P_2,\qquad \phi(1) = \begin{pmatrix}
1\\ &\ddots \\ &&1\\ &&&p
                                      \end{pmatrix},
\end{equation}
for $x\in L(H)$. Then 
\begin{equation}\label{SWWeq7.3}
Vx = Vfx = Vxf = VxV^*V = \phi(x)P_2V = \phi(x)V,\qquad x\in L(H),
\end{equation}
so
\begin{equation}\label{SWWeq7.4}
xV^* = V^*\phi(x),\qquad x\in L(H),
\end{equation}
by taking adjoints in \eqref{SWWeq7.3}. Note that \eqref{SWWeq7.4} is an equality of $1\times n$ row operators with entries from $\langle L(G), e_{L(H)}\rangle$. Let $\eta_j = v^*_j\xi\in \ell^2(G)$, $1\le j\le n$. For a fixed $j$,  the Kaplansky density theorem  allows us to choose a uniformly bounded net $(w_\alpha e_{L(H)})$ converging $*$-strongly to $v^*_je_{L(H)}$, where $w_\alpha\in L(G)$. For each $y\in L(H)$,
\begin{equation}\label{SWWeq7.5}
v^*_jy\xi = \lim_\alpha w_\alpha y\xi = \lim_\alpha Jy^*Jw_\alpha \xi = \eta_jy,
\end{equation}
and this equality then holds for each $j$, $1\le j\le n$.
Now apply \eqref{SWWeq7.4} to column vectors whose only non-zero entries are $\xi$ in the $j^{\rm th}$ component, $1\le j \le n$, and use \eqref{SWWeq7.5} to conclude that
\begin{equation}\label{SWWeq7.6}
x(\eta_1,\eta_2,\ldots,\eta_n) = (\eta_1,\eta_2,\ldots,\eta_n) \phi(x),\qquad x\in L(H),
\end{equation}
where the right action of $L(H)$ on $\ell^2(G)$ is used to define the multiplication on the right hand side of this equation. By putting $x=1$ in \eqref{SWWeq7.6}, we see that $\eta_n = \eta_np$, so \eqref{SWWeq7.6} can also be written as
\begin{equation}\label{SWWeq7.7}
u(\eta_1,\ldots, \eta_n) \phi(u^*) = (\eta_1,\ldots, \eta_n),\qquad u\in {\cl U}(L(H)). 
\end{equation}
Choose a sequence $(y_{1,m}\xi,\ldots,y_{n,m}\xi)$, $m\ge 1$, converging to $(\eta_1,\ldots,\eta_n)$ in $\|\cdot\|_2$-norm where $y_{i,m}\in L(G)$. The convex sets
\[
K_m = \ovl{\rm conv}^w\{u(y_{1,m},\ldots, y_{n,m})\phi(u^*)\colon\ u\in {\cl U}(L(H))\}
\]
are weakly  compact in $L(G)\times \cdots \times L(G)$, so the image in $\ell^2(G) \oplus\cdots\oplus \ell^2(G)$ is also weakly compact and weakly closed. Since $K_m$ is invariant for the action $u\cdot \phi(u^*)$, the unique element $(w_{1,m}\xi, \ldots, w_{n,m}\xi)\in K_m$ of minimal $\|\cdot\|_2$-norm, with $w_{i,m}\in L(G)$, satisfies
\begin{equation}\label{SWWeq7.8}
x(w_{1,m},\ldots, w_{n,m}) = (w_{1,m},\ldots, w_{n,m}) \phi(x),\qquad x\in L(H).
\end{equation}
Moreover, $\lim\limits_{m\to\infty} \|\eta_i-w_{i,m}\xi\|_2 = 0$ for $1\le i\le n$. It follows from \eqref{SWWeq7.8} that $\sum\limits^n_ {i=1} w_{i,m}w^*_{i,m}\in L(H)'\cap L(G) = {\bb C}1$ for each $m\ge 1$. Thus the $\|\cdot\|_2$-norm and operator norm agree for $(w_{1,m},\ldots, w_{n,m})$, $m\ge 1$. Since these converge to $(\eta_1,\ldots, \eta_n)$, they are bounded in $\|\cdot\|_2$-norm and hence in operator norm. By dropping to a subnet, we may further assume that they converge weakly to a row operator $(w_1,\ldots, w_n)\in L(G) \times\cdots\times L(G)$, whereupon $\eta_i = w_i\xi$ for $1\le i\le n$. From \eqref{SWWeq7.5}, we conclude that
\begin{equation}\label{SWWeq7.9}
v^*_jy\xi = \eta_jy = w_j\xi y = w_jy\xi,\qquad y\in L(H), 
\end{equation}
and so $v^*_je_{L(H)} = w_je_{L(H)}$ for $1\le j\le n$. Moreover, \eqref{SWWeq7.6} becomes
\begin{equation}\label{SWWeq7.10}
x(w_1,\ldots, w_n) = (w_1,\ldots, w_n)\phi(x),\qquad x\in L(H),
\end{equation}
and $w_np = w_n$, by putting $x=1$.

Let $\{g_iH\colon \ i\ge 1\}$ be a listing of the left $H$-cosets in $G$. Then there exist row operators $(z_{1,j},\ldots, z_{n,j})$, $j\ge 1$, with $z_{i,j} \in L(H)$ such that
\begin{equation}\label{SWWeq7.11}
(w_1,\ldots, w_n) = \sum^\infty_{j=1} \lambda_{g_j}(z_{1,j},\ldots, z_{n,j}),
\end{equation}
where the sum, which could be finite, converges in $\|\cdot\|_2$-norm, and $z_{n,j}p = z_{n,j}$. For each $h\in H$, \eqref{SWWeq7.10} gives
\begin{equation}\label{SWWeq7.12}
\lambda_h \sum^\infty_{j=1} \lambda_{g_j}(z_{1,j},\ldots, z_{n,j}) = \sum^\infty_{j=1} \lambda_{g_j}(z_{1,j},\ldots, z_{n,j}) \phi(\lambda_h).
\end{equation}

For convenience, write $Z_j = (z_{1,j},\ldots,z_{n,j})$ and suppose that the numbering has been chosen so that $\|Z_j\|_2\ge \|Z_{j+1}\|_2$, for $j\ge 1$, possible because $\|Z_j\|_2 \to 0$ as  $j\to\infty$. If $\|Z_j\|_2\ne 0$, then $S_j = \{i\colon \ \|Z_i\|_2 = \|Z_j\|_2\}$ is a finite set and, as in the proof of Theorem \ref{SWWthm6.2}, each $h\in H$ permutes the cosets $\{g_iH\colon \ i\in S_j\}$. We will now show that the number of non-zero $Z_j$'s must be at least $n$.

The range of $f$ is the range of $V^* = V^*e_{L(H)}$ and this operator is $(w_1e_{L(H)},\ldots, w_ne_{L(H)})$. Thus the range of $f$ is contained in the closure of
\[
\left\{\sum^n_{i=1} w_i\zeta_i\colon \ \zeta_i\in \ell^2(H)\right\}.
\]
Indeed, equality must hold since the projection onto this subspace is 
$$
\sum\limits^n_{i=1} w_ie_{L(H)}w^*_i = V^*V= f.
$$
If $Z_j$, $1\le j\le r<n$, are the only non-zero $Z_j$'s, then \eqref{SWWeq7.11} shows that the range of $f$ is contained in
\[
\left\{\sum^r_{i=1} \lambda_{g_i}\zeta_i\colon \ \zeta_i\in \ell^2(H)\right\}
\]
and the projection onto this space is $\sum\limits^r_{i=1} \lambda_{g_i} e_{L(H)} \lambda^{*}_{g_i}$ which has trace $r\le n-1$, contradicting $\text{Tr}(f)>n-1$. Thus we may pick an integer $N\ge n$ such that $\|Z_N\|_2 > \|Z_{N+1}\|_2$. Then each $h\in H$ permutes the left cosets $\{g_iH\colon \ 1\le i \le N\}$, so as in Theorem \ref{SWWthm6.2} there is a finite index normal subgroup $K$ of $H$ such that $g^{-1}_i\cdot g_i$ induces a homomorphism $\phi_i\colon \ K\to H$ for $1\le i\le N$. Since each $Z_j\ne 0$ for $1\le j\le N$, from \eqref{SWWeq7.11} we can find vectors $(\zeta_{1,j},\ldots, \zeta_{n,j})$, $1\le j\le N$, $\zeta_{i,j}\in \ell^2(H)$, such that $\sum\limits^n_{i=1} w_i\zeta_{i,j}$ has a non-zero $\lambda_{g_j}$-coefficient. A suitable linear combination then gives a vector $\sum\limits^\infty_{j=1} \lambda_{g_j}\zeta_j\in \text{Ran } f$, where $\zeta_j\in \ell^2(H)$, $j\ge 1$, and are non-zero for $1\le j \le N$. Pre-multiplication by $K$ and post-multiplication by $H$ allow us to find vectors in the range of $f$ whose first $N$ components are
\[
\sum^N_{i=1} \lambda_{g_i}\phi_i(\lambda_k) J\lambda^*_hJ\zeta_i,
\]
which we write in matrix form as
\begin{equation}\label{SWWeq7.13}
(\lambda_{g_1},\ldots, \lambda_{g_N}) \begin{pmatrix}
\phi_1(\lambda_k)J\lambda^*_hJ\\ &\ddots\\ &&\phi_N(\lambda_k)J\lambda^*_hJ 
                                      \end{pmatrix}
\begin{pmatrix}
\zeta_1\\ \vdots\\ \zeta_N
\end{pmatrix}.
\end{equation}

For $i\neq j$, $\phi_i\phi_j^{-1}(x)=\lambda_{g_i^{-1}g_j}x\lambda_{g_j^{-1}g_i}$ for $x\in L(G)$.  If there is a unitary $u\in L(H)$ with $\phi_i\phi_j^{-1}(y)=uyu^*$ for all $y\in L(K)$, then $u^*\lambda_{g_i^{-1}g_j}\in L(K)'\cap L(G)=\mathbb C1$. Hence $g_i^{-1}g_j\in H$, which is a contradiction as $g_iH\neq g_jH$.  We can now apply Lemma \ref{SWWlem2.2} to deduce that the diagonal matrices in (\ref{SWWeq7.13}) generate the von Neumann algebra
\[
\left\{\begin{pmatrix}
t_1\\ &\ddots\\ &&t_N 
       \end{pmatrix}\colon \ t_i\in  B(\ell^2(H))\right\}.
\]
Since $\zeta_i\ne 0$ for $1\le i\le N$, we see that
\[
\text{span}\left\{\begin{pmatrix}
\phi_1(\lambda_k)J\lambda^*_hJ\\ &\ddots\\ &&\phi_N(\lambda_k)J\lambda^*_hJ
                  \end{pmatrix}
\begin{pmatrix}
\zeta_1\\ \vdots\\ \zeta_N
\end{pmatrix}\colon \ k\in K,\quad h\in H\right\}
\]
is dense in the direct sum of $N$ copies of $\ell^2(H)$. If $\tilde f = \sum\limits^N_{i=1} \lambda_{g_i} e_{L(H)}\lambda^*_{g_i}$, then $\tilde f$ is a projection of trace $N$. The range projection of $\tilde ff$ is $\tilde f$ while the range projection of $f\tilde f$ lies under $f$. Since these range projections are equivalent in $\langle L(G), e_{L(H)}\rangle$, we conclude that $\text{Tr}(\tilde f)\le \text{Tr}(f)$. Thus
\begin{equation}\label{SWWeq7.14}
n\le N = \text{Tr}(\tilde f) \le \text{Tr}(f) = n-1+\mu,
\end{equation}
forcing $\mu=1$, and $N=n$. In particular, no choice of $N>n$ was possible. Thus only $\lambda_{g_j}$ terms for $j\le n$ appear in \eqref{SWWeq7.11} and so $f\le \tilde f$. Equality of the traces then gives $f=\tilde f$, and the result follows.
\end{proof}

\begin{cor}\label{SWWcor7.2}
Let $H\subseteq G$ be an inclusion of countable discrete groups such that $L(H)' \cap L(G) = {\bb C}1$. Then the projections in $L(H)' \cap \langle L(G), e_{L(H)}\rangle$ of finite trace generate an abelian algebra.
\end{cor}

\begin{proof}
Let $\{g_iH\colon\ i\ge 1\}$ be a listing of distinct left $H$-cosets in $G$. Then the projections $\{\lambda_{g_i}e_{L(H)}\lambda^*_{g_i}\colon \ i\ge 1\}$ are pairwise orthogonal, and suitable finite sums of these account for all finite trace projections in $L(H)' \cap \langle L(G), e_{L(H)}\rangle$, by Theorem \ref{SWWthm7.1}. The result is immediate from this. 
\end{proof}

We note that Corollary \ref{SWWcor7.2} is special to the group-subgroup situation. If a finite group $G$ acts on a factor $M$ and we let $N$ be $M^G$, then we have already seen that $N'\cap  \langle M,e_N\rangle$ is isomorphic to the group algebra ${\bb C}G$ which can be highly non-abelian. Even for $L(H)\subseteq L(G)$, the situation for infinite trace projections can be very complicated, as the next example shows.

\begin{exm}\label{SWWexm7.3}
Let $G = {\bb F}_3$, generated by $\{a,b,c\}$, and let $H = {\bb F}_2$, generated by $\{a,b\}$. We will show that $L(H)' \cap \langle L(G), e_{L(H)}\rangle$ is the direct sum of ${\bb C}$ and a $\text{II}_\infty$ factor. First consider two reduced words $w_1,w_2\in {\bb F}_3$ which begin and end with non-zero powers of $c$. For $i=1,2$, let $f_i$ be the projection onto the $L(H)$-bimodule generated by $w_i$ in $\ell^2(G)$. Then $f_1$ and $f_2$ are in $L(H)'\cap \langle L(G),e_{L(H)}\rangle$ and we establish that they are equivalent in this algebra. There can be no cancellation of letters in the words $x_iw_1y_i$ when $x_i,y_i\in {\bb F}_2$, $1\leq i\leq n$,    are reduced words. Then $\{x_iw_1y_i\colon 1\le i\le n\}$ is an orthogonal set so, for scalars $\alpha_i$,
\begin{align}\label{SWWeq7.15}
\left\|\sum^n_{i=1} \alpha_ix_iw_1y_i\right\|^2_2 = \sum^n_{i=1} |\alpha_i|^2, 
\end{align}
with a similar equality for $w_2$ replacing $w_1$. Thus there is an isometric map $v\colon \ \text{Ran } f_1\to \text{Ran } f_2$ defined on a generating set of vectors by 
\begin{equation}\label{SWWeqw7.16}
v(xw_1y) = xw_2y,\qquad x,y\in {\bb F}_2,
\end{equation}
and this extends to $\ell^2(G)$ by setting $v=0$ on $(\text{Ran } f_1)^\bot$. The definition makes it clear that $v$ commutes with the left and right actions of $L(H)$, so $v\in L(H)' \cap \langle L(G), e_{L(H)}\rangle$, and $v^*v =f_1$, $vv^*=f_2$. Thus $1-e_{L(H)}$ is an orthogonal sum of equivalent projections corresponding to words beginning and ending with a non-zero power of $c$. 
Since there are infinitely many such words, this will contribute a tensor factor of 
$B(\ell^2(\bb N))$.
Thus it suffices to consider projections in $L(H)'\cap \langle L(G), e_{L(H)}\rangle$ whose ranges are contained in the $L(H)$-bimodule $X_c$ generated by $c$. We denote the projection onto $X_c$ by $f_c$.

Consider the $\text{II}_1$ factor $L(H)\ \ovltimes\ JL(H)J$ acting in standard form on $L^2(L(H)\ \ovltimes\ JL(H)J)$. If we define an operator ${U}\colon \ L^2(L(H)\ \ovltimes\ JL(H)J)\to X_c$ by
\begin{equation}\label{SWWeq7.17}
{U}\left(\sum^n_{i=1} \alpha_i(x_i \otimes Jy_iJ)\right) = \sum^n_{i=1} \alpha_ix_icy^{-1}_i
\end{equation}
for words $x_i,y_i\in {\bb F}_2$ and scalars $\alpha_i$, then the calculation of \eqref{SWWeq7.15} shows that ${U}$ extends to a well defined unitary between these Hilbert spaces. Since 
\begin{equation}\label{SWWeq7.18}
{U}(x \otimes JyJ){U}^*(c) =  xcy^{-1}
\end{equation}
for words $x,y\in {\bb F}_2$, we see that ${U}(L(H)\ \ovltimes\ JL(H)J)U^*$ is the von
Neumann algebra
 of operators on $X_c$ generated by left and right multiplications by elements of $L(H)$. Thus ${U}$ conjugates $(L(H)\ \ovltimes\ JL(H)J)'$ to $f_c(L(H)'\cap \langle L(G), e_{L(H)}\rangle)f_c$, from which it follows that
\begin{align}
L(H)' \cap \langle L(G),e_{L(H)}\rangle &\cong {\bb C}e_{L(H)} \oplus (1-e_{L(H)})(L(H)'\cap \langle L(G), e_{L(H)}\rangle)\notag\\
&\cong {\bb C}e_{L(H)} \oplus (L(H)\ \ovltimes\ JL(H) J)'\ \ovltimes\ B(\ell^2({\bb N}))\notag\\
\label{SWWeq7.19}
&\cong {\bb C}e_{L(H)} \oplus (JL(H)J)\ \ovltimes\ L(H)\ \ovltimes\ B(\ell^2({\bb N})),
\end{align}
and the latter tensor product algebra is a $\text{II}_\infty$ factor. We also conclude that there are no minimal projections in $L(H)'\cap \langle L(G), e_{L(H)}\rangle$ under $1-e_{L(H)}$.
$\hfill\square$
\end{exm}

Example \ref{SWWexm7.3} showed that $L(H)'\cap \langle L(G),e_H\rangle$ could contain a $\text{II}_\infty$ factor. Our final example shows that, even in the singular infinite index case, this algebra can also be atomic, abelian and generated by its minimal projections of finite trace. Furthermore, the traces of these minimal projections can be uniformly bounded. We note that any countable discrete group $G$ can act on ${\bb F}_{|G|}$ by outer automorphisms. Index the generators of ${\bb F}_{|G|}$ by $\{g_t\colon \ t\in G\}$ and let $\beta_s\in \text{Aut}({\bb F}_{|G|})$ be defined on generators by $g_t\mapsto g_{st}$, $s,t\in G$. The semidirect product $\mathbb F_{|G|}\rtimes_\beta G$ is a countable I.C.C. group.

\begin{exm}\label{SWWexm7.4}
Let ${\bb Z}_2$ act on ${\bb Z}$ by
\begin{equation}\label{SWWeq7.20}
\alpha_m(n) = (-1)^mn,\qquad n\in {\bb Z},\quad m\in {\bb Z}_2, 
\end{equation}
and then let ${\bb Z}\rtimes_\alpha {\bb Z}_2$ act on ${\bb F}_\infty$ by an action $\beta$ as described above. Set $G = {\bb F}_\infty \rtimes_\beta ({\bb Z}\rtimes_\alpha {\bb Z}_2)$ and let $H$ be the subgroup generated by ${\bb F}_\infty$ and ${\bb Z}_2$. Each $g\in G$ has infinitely many $H$-conjugates and so $L(H)$ is irreducible in $L(G)$. Any $g\in G\backslash H$ contains a non-zero group element $n\in {\bb Z}$, and then properties of the semidirect product show that the double coset $HgH$ is a union of two left cosets generated by $\pm n\in {\bb Z}$. By Theorem \ref{SWWthm7.1}, we see that each of these double cosets corresponds to a minimal projection in $L(H)' \cap \langle L(G),e_{L(H)}\rangle$ of trace 2, so this algebra is abelian and any projection in it under $1-e_{L(H)}$ is an orthogonal sum of projections of trace $2$.

Many variations on this theme are possible. Replace ${\bb Z}_2$ by a group of order $n$ and replace ${\bb Z}$ by an infinite group on which it acts. The minimal projections will then all have integer trace bounded by $n$.
\end{exm}


\begin{thebibliography}{10}

\bibitem{Chifan.Normalisers}
I.~Chifan.
\newblock On the normalizing algebra of a masa in a $\mathrm{II}_1$ factor.
\newblock Preprint, arXiv:math.OA/0606225, 2006.

\bibitem{Christensen.Subalgebras}
E.~Christensen.
\newblock Subalgebras of a finite algebra.
\newblock {\em Math. Ann.}, 243(1):17--29, 1979.

\bibitem{Dixmier.Masa}
J.~Dixmier.
\newblock Sous-anneaux ab\'eliens maximaux dans les facteurs de type fini.
\newblock {\em Ann. of Math. (2)}, 59:279--286, 1954.

\bibitem{Sinclair.FreePuk}
K.~J. Dykema, A.~M. Sinclair, and R.~R. Smith.
\newblock Values of the {Puk\'anszky} invariant in free group factors and the
  hyperfinite factor.
\newblock {\em J. Funct. Anal.}, 240:373--398, 2006.

\bibitem{Popa.IPP}
A.~Ioana, J.~Peterson, and S.~Popa.
\newblock Amalgamated free products of $w$-rigid factors and calculation of
  their symmetry groups.
\newblock Acta Math., to appear. arXiv.OA/0505589.

\bibitem{Jones.RegularSubfactors}
V.~Jones.
\newblock Sur la conjugaison de sous-facteurs de facteurs de type ${\rm{II}}_1$.
\newblock {\em C. R. Acad. Sci. Paris S\'er. A-B}, 284(11):A597--A598, 1977.

\bibitem{Jones.PropertiesMasas}
V.~Jones and S.~Popa.
\newblock Some properties of {MASA}s in factors.
\newblock In {\em Invariant subspaces and other topics
  (Timi{\c{s}}oara/Herculane, 1981)}, volume~6 of {\em Operator Theory: Adv.
  Appl.}, pages 89--102. Birkh\"auser, Basel, 1982.

\bibitem{Jones.SubfactorsBook}
V.~Jones and V.~S. Sunder.
\newblock {\em Introduction to subfactors}, volume 234 of {\em London
  Mathematical Society Lecture Note Series}.
\newblock Cambridge University Press, Cambridge, 1997.

\bibitem{Jones.Index}
V.~F.~R. Jones.
\newblock Index for subfactors.
\newblock {\em Invent. Math.}, 72(1):1--25, 1983.

\bibitem{Kosaki.Bimodules}
H.~Kosaki and S.~Yamagami.
\newblock Irreducible bimodules associated with crossed product algebras.
\newblock {\em Internat. J. Math.}, 3(5):661--676, 1992.

\bibitem{MvN.4}
F.~J. Murray and J.~von Neumann.
\newblock On rings of operators. {IV}.
\newblock {\em Ann. of Math. (2)}, 44:716--808, 1943.

\bibitem{Nakamura.GaloisFinite}
M.~Nakamura and Z.~Takeda.
\newblock A {G}alois theory for finite factors.
\newblock {\em Proc. Japan Acad.}, 36:258--260, 1960.

\bibitem{Nakamura.Galois}
M.~Nakamura and Z.~Takeda.
\newblock On the fundamental theorem of the {G}alois theory for finite factors.
\newblock {\em Proc. Japan Acad.}, 36:313--318, 1960.

\bibitem{Popa.Entropy}
M.~Pimsner and S.~Popa.
\newblock Entropy and index for subfactors.
\newblock {\em Ann. Sci. \'Ecole Norm. Sup. (4)}, 19(1):57--106, 1986.

\bibitem{Popa.NotesCartan}
S.~Popa.
\newblock Notes on {C}artan subalgebras in type {${\rm II}\sb 1$} factors.
\newblock {\em Math. Scand.}, 57(1):171--188, 1985.

\bibitem{Popa.CBMSNotes}
S.~Popa.
\newblock {\em Classification of subfactors and their endomorphisms}, volume~86
  of {\em CBMS Regional Conference Series in Mathematics}.
\newblock Published for the Conference Board of the Mathematical Sciences,
  Washington, DC, 1995.

\bibitem{Popa.Betti}
S.~Popa.
\newblock On a class of type {${\rm II}\sb 1$} factors with {B}etti numbers
  invariants.
\newblock {\em Ann. of Math. (2)}, 163(3):809--899, 2006.

\bibitem{Popa.StrongRigidity1}
S.~Popa.
\newblock Strong rigidity of {$\rm{II}_1$} factors arising from malleable
  actions of {$w$}-rigid groups, I.
\newblock {\em Invent. Math.}, 165:369---408, 2006.

\bibitem{Popa.StrongRigidity2}
S.~Popa.
\newblock Strong rigidity of {$\rm{II}_1$} factors arising from malleable
  actions of {$w$}-rigid groups, II.
\newblock {\em Invent. Math.}, 165:409---451, 2006.

\bibitem{Sinclair.PertSubalg}
S.~Popa, A.~M. Sinclair, and R.~R. Smith.
\newblock Perturbations of subalgebras of type {II{$\sb 1$}} factors.
\newblock {\em J. Funct. Anal.}, 213(2):346--379, 2004.

\bibitem{Sinclair.StrongSing2}
G.~Robertson, A.~M. Sinclair, and R.~R. Smith.
\newblock Strong singularity for subalgebras of finite factors.
\newblock {\em Internat. J. Math.}, 14(3):235--258, 2003.

\bibitem{Sinclair.MasaBook}
A.~M. Sinclair and R.~R. Smith.
\newblock Finite von {N}eumann algbras and masas.
\newblock  Book manuscript, Cambridge University Press, to appear.

\bibitem{Sinclair.StrongSing}
A.~M. Sinclair and R.~R. Smith.
\newblock Strongly singular masas in type {$\rm II\sb 1$} factors.
\newblock {\em Geom. Funct. Anal.}, 12(1):199--216, 2002.

\bibitem{Sinclair.Puk}
A.~M. Sinclair and R.~R. Smith.
\newblock The {P}uk\'anszky invariant for masas in group von {N}eumann factors.
\newblock {\em Illinois J. Math.}, 49(2):325--343 (electronic), 2005.

\bibitem{Saw.StrongSing}
A.~M. Sinclair, R.~R. Smith, S.~A. White, and A.~Wiggins.
\newblock Strong singularity of singular masas. Illinois J. Math., to appear.
\newblock arXiv:math.OA/0601594.

\bibitem{Saw.Length}
S.~A. White and A.~Wiggins.
\newblock Semi-regular masas of transfinite length.
\newblock Internat. J. Math., to appear. arXiv:math.OA/0611615.

\bibitem{Wiggins.StrongSing}
A.~Wiggins.
\newblock Strong singularity for subfactors of a $\mathrm{II}_1$ factor.
\newblock Preprint, arXiv:math.OA/0703673, 2007.

\end{thebibliography}
\end{document}